\documentclass[11pt]{article}
\usepackage{amssymb}
\usepackage{amsmath}
\usepackage{amsopn}
\usepackage{pifont}
\usepackage{color}
\usepackage{mathrsfs}
\usepackage{amsthm}
\usepackage{ulem}

\numberwithin{equation}{section}

\theoremstyle{plain}
\newtheorem{thm}{Theorem}[section]
\newtheorem{prop}{Proposition}[section]
\newtheorem{lemma}{Lemma}[section]

\theoremstyle{definition}

\newtheorem{remark}{Remark}[section]

\def\<{\langle}
\def\>{\rangle}

\def\CC{{\mathbb C}}
\def\ZZ{{\mathbb Z}}

\def\r#1{(\ref{#1})}
\def\ot{\otimes}
\def\sk#1{\left(#1\right)}

\def\Pep{{P}^+_e}
\def\Pem{{P}^-_e}
\def\Pepm{{P}^\pm_e}

\def\Pfpm{{P}_f^\pm}
\def\Pfp{{P}^+_f}
\def\Pfm{{P}^-_f}
\def\RR{{\rm R}}
\def\LL{{\rm L}}
\def\hLL{\hat{\rm L}}

\def\E{{\sf e}}

\def\proof{\noindent{\it Proof.}\ }

\def\gle{g}
\def\gri{\tilde g}

\def\qed{\hfill\nobreak\hbox{$\square$}\par\medbreak}
\def\Id{\mathbb{I}}
\def\Per{\mathbb{P}}
\def\Qer{\mathbb{Q}}
\def\Der{\mathbb{D}}
\def\Derd{\tilde{\mathbb{D}}}
\def\Rer{\mathbb{R}}

\def\gaf{\tilde{\mathfrak{g}}}
\def\gaff{\widehat{\mathfrak{g}}}

\def\qf{{\sf q}}
\def\pf{{\sf p}}
\def\Dn#1{D^{(2)}_{#1}}
\def\Ut{U}
\def\FF{\mathrm{F}}
\def\EE{\mathrm{E}}
\def\tFF{\tilde{\mathrm{F}}}
\def\tEE{\tilde{\mathrm{E}}}

\def\TT{\mathbb{U}}
\def\Fn{\mathcal{F}}
\def\En{\mathcal{E}}
\def\PPd{\mathcal{P}}
\def\QQd{\mathcal{Q}}

\def\MM{\mathrm{M}}
\def\Lf{\mathbb{L}}

\def\QQdd{{\sf Q}_2}

\def\adm{{\sf a}}
\def\bdm{{\sf b}}
\def\Gap{\mathfrak{q}}
\def\QQff{{\sf Q}_1}
\def\fgo{\mathfrak{f}}
\def\ggo{\mathfrak{g}}

\begin{document}

\thispagestyle{empty}
\setcounter{page}{0}

\vspace{12pt}

\begin{center}
\begin{Large}
{\bf On the  R-matrix realization of the quantum\\[2mm]
loop algebra. The case of $U_q(D^{(2)}_n)$}
\end{Large}

\vspace{20pt}

\begin{large}
A.~Liashyk${}^{a}$ and  S.~Pakuliak${}^{b}$
\end{large}

\vspace{10mm}

${}^a$ {\it Beijing Institute of Mathematical Sciences and Applications (BIMSA),\\
Beijing, China\\
E-mail: liashyk@bimsa.cn}

\vspace{2mm}

${}^b$ {\it Laboratoire d'Annecy-le-Vieux de Physique Théorique (LAPTh)\\
Chemin de Bellevue, BP 110, F-74941, Annecy-le-Vieux cedex, France\\
E-mail: pakuliak@lapth.cnrs.fr}

\end{center}

\noindent{\footnotesize {\bf Abstract.}
The connection between the R-matrix realization and Drinfeld's realization of the quantum loop algebra $U_q(D^{(2)}_n)$ is considered using the Gaussian decomposition approach proposed by J. Ding and I. B. Frenkel.
Our main result is a description of the embedding $U_q(\Dn{n-1})\hookrightarrow  U_q(\Dn{n})$ that underlies this connection.
Explicit relations between all Gaussian coordinates of the $\LL$-operators and the currents are presented.
}

\tableofcontents

\section{Introduction}

It is known  that quantum affine algebras and Yangians \cite{Drinfeld1988quantumgroups} have three equivalent
realizations.
The first realization is given in terms of the finite number of
Chevalley generators \cite{Drinfeld1988quantumgroups, chari1995guide}.
The second realization, which replicates  the Quantum Inverse Scattering Method
\cite{SklyaninTakhtadzhyanFaddeev1979QISM} is related
to the description in terms of $\LL$-operators and $\RR$-matrices
   \cite{ReshetikhinSklyanin1990quantumcurrentgroup}.
 The third realization was formulated in
\cite{Drinfeld1988NewRealisation} and give a description of the
Yangians and the quantum affine algebras in terms of formal generating
series or {\it currents}.

A successful approach to prove an equivalence of second and third realizations
for the quantum affine algebra $U_q(\widehat{\mathfrak{gl}}(n))$ was
proposed in \cite{DingFrenkel1993isomorphism}. The main tool
 was so-called {\it Gaussian coordinates}  of the
fundamental $\LL$-operators. It was shown in \cite{DingFrenkel1993isomorphism}
that the linear combinations of the simple root Gaussian coordinates coincide
with Drinfeld's generating series (currents) of the third realization.
The reason for this success is understood today because of the known structure
of the universal $\mathcal{R}$-matrices  for quantum affine algebras and
Yangian doubles \cite{KhoroshkinTolstoy1991UniversalRmatrix}
written in terms of Cartan-Weyl generators.
$\LL$-operators used in the second realization can be
obtained by evaluating universal $\mathcal{R}$-matrix in the fundamental
representation in one of the tensor components.
It was shown in \cite{KHOROSHKIN1993445}
that Cartan-Weyl basis for quantum affine algebras coincides with
Drinfeld's current generators and this explains the main result of
the paper \cite{DingFrenkel1993isomorphism}.  This is a reason
why we call the last or Drinfeld 'new' realization as {\it Cartan-Weyl (CW) realization}.

Until 2017, most papers devoted to isomorphisms of different realizations of the Yangians and the quantum affine algebras were related to the $A$-series algebras and their super-symmetric generalizations.

In 2017, Gaussian decomposition of the fundamental $\LL$-operators was used
to describe isomorphisms between second and third realizations
of the Yangian doubles of $B$, $C$ and $D$ series algebras \cite{JingLiuMolev2018Y(BCD)}.
Earlier the similar approach was used in  \cite{JingLiu2013Y(so(3))} for the
Yangian $Y(\mathfrak{so}(3))$ associated with
simplest algebra $\mathfrak{so}(3)$ of the $B$-series.

Isomorphism between $\LL$-operator and CW descriptions of the
quantum affine algebras were presented in the papers
\cite{JingLiuMolev2020Uq(C),JingLiuMolev2020Uq(BD)} and extension of these
results to the twisted $A$-series algebras can be found in
\cite{LiashykPakuliak2022U_q(ABCD)}. The case of
$U_q(A^{(2)}_2)$ was considered in \cite{Shapiro2010three}.

Construction proposed in the paper \cite{DingFrenkel1993isomorphism} describes the connection between Gaussian coordinates and currents associated only with simple roots.
The conception of ``composed current'' related to non-simple roots was introduced in \cite{DingKhoroshkin2000Weylgroup}.
Connection of the ``composed currents'' to all the Gaussian coordinates was explained in \cite{KyotoPaper} for the case of $U_q(\widehat{\mathfrak{gl}}(n))$. For other series, this connection was given in the paper \cite{LiashykPakuliak2020Y(BCD)} for Yangians and in \cite{LiashykPakuliak2022U_q(ABCD)} for quantum affine algebras.

A classification of the solutions to the classical
and quantum Yang-Baxter equations was described for the case of nonexceptional quantum
affine Lie algebras in \cite{Jimbo1986quantum}.
The papers \cite{JingLiuMolev2020Uq(C), JingLiuMolev2020Uq(BD), LiashykPakuliak2022U_q(ABCD)} describe the connection between the different realizations for all the cases listed in \cite{Jimbo1986quantum} except the case $\Dn{n}$. In this paper, we fill this gap.

Our main tool is the embedding $U_q(\Dn{n-1})\hookrightarrow  U_q(\Dn{n})$ (see theorem \ref{main-st}).
We use this theorem to investigate the connection between the second and the third realizations of quantum loop algebra $U_q(\Dn{n})$.
To define quantum loop algebra $U_q(\Dn{n})$, we use $\RR$-matrix \eqref{tRalt} and then apply a Ding-Frenkel construction.

The paper is composed as follows. In section~\ref{R-matrix} we describe
$\RR$-matrix \eqref{tRalt}, its relation to $\RR^{\rm J}(u, v)$ \eqref{RJmat} found by M. Jimbo, and its properties which are proved in appendix~\ref{sec: Crossing symmetries and pole structure}.
$\RR$-matrix \eqref{tRalt} is used to define the quantum loop algebra $U_q(\Dn{n})$ in section~\ref{Algebra}.
We understand this algebra as a collection of matrix elements of the fundamental $\LL$-operators
which are formal series with respect to the spectral parameters $u$ and $u^{-1}$.
We discuss the center of the algebra in  section~\ref{sec-cent}.
The Gaussian coordinates of the fundamental $\LL$-operators are introduced in section~\ref{GCsect}.
Section~\ref{Algebra} is concluded by discussion of the normal ordering
of the Gaussian coordinates as a generating series of the
CW generators and introduces  projections on the
intersections of the different types of the Borel subalgebras
 \cite{EnriquezKhoroshkinPakuliak2007Projections}.
Section~\ref{embed} is devoted to a description of embedding
of the algebra $U_q(\Dn{n-1})$ in $U_q(\Dn{n})$.
This embedding is formulated as theorem~\ref{main-st}, and the proof
of this theorem is based on embedding relations for $\RR$-matrix \r{tRalt}
formulated in this section and proved in appendix~\ref{ApD}.
Section~\ref{cur-real} describes the realizations of
quantum loop algebra $U_q (D^{(2)_n})$ in terms of the currents. Small rank algebras with $n=2,3$ are considered in more detail.
Relations between all Gaussian coordinates and the currents are given in this section. Appendix~\ref{ApA}
recalls Jimbo's $\RR$-matrix \cite{Jimbo1986quantum}.
Serre relations are interpreted in appendix~\ref{ApE} as the commutation relations between composed currents.

\section{$\RR$-matrix}\label{R-matrix}

Let $\mathfrak{g}$ be  the  Lie algebra  $\mathfrak{o}_{2n}$.
Let $\gaff$  be twisted affine algebra $D^{(2)}_{n}$.
By  $\gaf$ we denote the loop algebra, which is obtained from
the affine algebra $\gaff$ setting to zero a central charge element.

Let $q\in \CC$ be an arbitrary complex number not equal to zero or root of unity.
We consider standard quantum deformation  $U_q(\gaf)$ of the universal enveloping
algebra  $U(\gaf)$ \cite{Drinfeld1988quantumgroups}.

Let $\E_{ij}$ be  $N\times N$ matrix units for $1\leq i,j\leq N = 2n$ and we denote
\begin{equation}\label{prime}
i'=N+1-i,\quad 1\leq i\leq N\,.
\end{equation}
To describe $U_q(\Dn{n})$-invariant  $\RR$-matrix
\cite{Jimbo1986quantum}
one needs a parameter $\xi$
\begin{equation}\label{kappa}
\xi = q^{-n+1}
\end{equation}
and  the map $\bar\imath$ for $i=1,\ldots,N$
\begin{equation}\label{map}
\bar\imath=
\begin{cases}
  n - \frac{1}{2} - i, &\text{ for } i < n,\\
  0, &\text{ for } i = n, n+1, \\
  n + \frac{3}{2} - i, &\text{ for } i > n+1.\\
\end{cases}
\end{equation}
For example, $\bar 1=n-3/2$, $\bar 2=n-5/2$, $\bar n=0$, etc.

Let us introduce  rational functions
depending on $q$ and   spectral parameters
$u$ and $v$
\begin{equation}\label{rat-fun}
f(u,v)=\frac{u q-vq^{-1}}{u-v}\,,\quad
\gle(u,v)=\frac{\gamma_q u}{u-v}\,,\quad
\gri(u,v)=\frac{\gamma_q v}{u-v},\quad \gamma_q=q-q^{-1}\,.
\end{equation}

For all $1\leq i,j\leq 2n$ we set
\begin{equation}\label{p-fun}
\pf_{ij}(u,v)=\begin{cases}f(u,v)-1,&i=j\,,\\  \gri(u,v),&i<j\,,\\ \gle(u,v),&i>j\,,
\end{cases}
\end{equation}
\begin{equation}\label{q-fun}
    \qf_{ij}(u,v|\xi) = q^{\bar{\imath} - \bar{\jmath}} \pf_{ji}(v\xi, u) -
    \alpha_q \delta_{\bar\imath, \bar{\imath'}}\delta_{ij}\,,
\end{equation}
where parameter $\alpha_q$ is
\begin{equation}\label{Tta}
\alpha_q=q-2+q^{-1}\,.
\end{equation}

Let $\Per(u,v)$ and $\Qer^{\rm J}(u,v|\xi)$ be operators in ${\rm End}(\CC^N\otimes \CC^N)$
\begin{equation}\label{PPuv}
\Per(u,v)= \sum_{1\leq i,j\leq N} \pf_{ij}(u,v)\ \E_{ji}\ot \E_{ij}\,,
\end{equation}
\begin{equation}\label{QQuv}
\Qer^{\rm J}(u,v|\xi)=
\sum_{1\leq i,j\leq N} \qf_{ij}(u,v|\xi)\  \E_{i'j'}\ot \E_{ij}\,.
\end{equation}
For any $X\in{\rm End}(\CC^N)$ a transposed matrix $X^{\rm t}$ is defined as
\begin{equation}\label{trans}
(X^{\rm t})_{ij}=X_{j'i'}\,.
\end{equation}
One can check that
$\Qer^{\rm J}(u,v)$ and $\Per(u,v)$ are connected by the relation
\begin{equation}\label{con1}
    \Qer^{\rm J}_{12}(u,v|\xi) = \Der_2\ \Per_{12}\ \Per_{12}(v\xi, u)^{{\rm t}_1}\
    \Per_{12}\ \Der^{-1}_2 - \alpha_q \sum_{i = n, n+1} \E_{ii}\ot\E_{ii}\,,
\end{equation}
where $\Der$ is a diagonal matrix
\begin{equation}\label{Der}
\Der_{ij}=\delta_{ij} q^{\bar\imath}
\end{equation}
and $\Der_2=\Id\ot\Der$. $\Per_{12}$ is a permutation operator
\begin{equation*}
\Per_{12}=\sum_{1\leq i,j\leq N}\E_{ij}\ot\E_{ji}
\end{equation*}
and one can verify that \r{con1} is equivalent to the relations \r{q-fun}
which express $\qf_{ij}(u,v|\xi)$ through $\pf_{ij}(u,v)$.

Let $\Id$ be identity operator in ${\rm End}(\CC^N)$ and matrix $\Ut\in{\rm End}(\CC^N)$ is
\begin{equation}\label{Ut}
\Ut_{ij}=\delta_{ij}\delta_{i,j\not=n,n+1}+\delta_{ij'}\delta_{i,j=n,n+1}\,.
\end{equation}
Here, $\delta_{\rm condition}$ is equal to 1 if the statement of the condition is true and 0 otherwise.

It is proved in appendix \ref{ApA} that Jimbo's $\RR$-matrix $\RR^{\rm J}(u,v)$ \cite{Jimbo1986quantum} has a
presentation
\begin{equation}\label{D2Rmat}
\begin{split}
\RR^{\rm J}(u,v)\ &=\frac12\Big(\Id \otimes \Id +\Per(u,v) +  \Qer^{\rm J}(u,v|\xi)\Big)+\\
&+ \frac12(\Ut\ot\Id)\cdot \Big(\Id \otimes \Id+\Per(-u,v) +  \Qer^{\rm J}(-u,v|\xi)\Big)
\cdot( \Ut\ot\Id)  \,,
 \end{split}
\end{equation}
where matrices $\Per(u,v)$, $\Qer^{\rm J}(u,v|\xi)$ and $\Ut$ are given by equations
\r{PPuv}, \r{QQuv} and \r{Ut} respectively.

In this work, we use the slightly different R-matrix to define loop algebra  $U_q(\Dn{n})$.
\begin{prop}\label{D-inv-R}
  $\RR$-matrix for $U_q(\Dn{n})$ has the form
  
\begin{multline}\label{tRalt}
  \RR(u,v) =(\TT\ot\TT)\cdot \RR^{\rm J}(u,v)\cdot (\TT^{-1}\ot\TT^{-1}) 
             = \\ \frac12\Big(\Id \otimes \Id +\Per(u,v) +  \Qer(u,v|\xi)\Big)  +\frac12(\Ut\ot\Id)\cdot \Big(\Id \otimes \Id+\Per(-u,v) +  \Qer(-u,v|\xi)\Big)\cdot( \Ut\ot\Id),
\end{multline}
where
\begin{equation}\label{Qeven}
\begin{aligned}
    \Qer_{12}(u,v|\xi) = (\TT\ot\TT) \Der_2\ \Per_{12}\ \Per_{12}(v\xi, u)^{{\rm t}_1}  \Per_{12}\ \Der^{-1}_2 (\TT^{-1}\ot\TT^{-1}) -  \alpha_q \sum_{i = n, n+1} \E_{ii}\ot\E_{ii} \, + X_{12},   
\end{aligned}
\end{equation}
with 
\begin{equation}\label{eq: X}
 X_{12} = \frac{\gamma_q}{8\xi^{1/2}}\Big(\E_{n,n+1} + \E_{n+1,n}  \Big) \ot
   \left((\xi-1) \Big(\E_{n,n} - \E_{n+1,n+1}\Big)   
   +(\xi^{1/2}-1)^2\Big( \E_{n+1,n} - \E_{n,n+1} \Big)   \right),
\end{equation}
and
\begin{equation}\label{Tm1}
\TT_{ij}=\delta_{ij}\delta_{i,j\not=n,n+1}+\big(a_+\delta_{ij}+a_-\delta_{ij'})\delta_{i,j=n,n+1}\,,
\quad a_\pm=\frac{1\pm\xi^{-1/2}}{2q^{1/4}}.
\end{equation}
\end{prop}

For $\PPd(u,v)$ defined in \eqref{PP14} we have
\begin{equation}
 ( \TT \otimes \TT )\, \PPd(u,v) \, ( \TT^{-1} \otimes \TT^{-1} ) = \PPd(u,v) + X,
\end{equation}
Note that $X$ satisfies the equality
\begin{equation}
  (U \otimes \Id)X(U \otimes \Id) = X.
\end{equation}

Using representation \eqref{con1} for $\Qer^{\rm J}$ we immediately obtain the first line in the equality \eqref{Qeven} as a consequence of conjugation by the $\TT \otimes \TT$.
The second line in \eqref{Qeven} is literally $\tfrac{1}{2} X$ \eqref{eq: X}.

Note that the constant term in the second line of the formula \eqref{Qeven} is fixed up to term $\chi$, which satisfies the equality
$\chi + (U \otimes \Id) \chi (U \otimes \Id) = 0$. \qed

Twisting of Jimbo's $R$-matrix  \eqref{D2Rmat} is necessary to have diagonal $\mathbf{K}$-matrix in the Gaussian decomposition \eqref{Gauss3}.
Original Jimbo's R-matrix \eqref{D2Rmat} is not compatible with the diagonal $\mathbf{K}$-matrix (see appendix \ref{sec:GD Jimbo R-matrix}).

Using $\RR$-matrix \r{tRalt} we can prove embedding of $U_q(\Dn{n-1})$ into $U_q(\Dn{n})$ presented in section~\ref{sec:embegging theorem}.
The meaning of the embedding theorem~\ref{main-st} is that in order to obtain
 the commutation relations between Gaussian coordinates of $\LL$-operators
 or between corresponding currents, it is sufficient to get these relations
 for the small rank algebras.  Also, this theorem implies that the structure
 of the internal part of $U_q(\Dn{n})$ $\LL$-operator should be the same
 as for the $\LL$-operator of the algebra $U_q(\Dn{n-1})$. This will be demonstrated
 in section~\ref{cur-real} explicitly for $n = 2,3$.

\subsection{Properties of R-matrix}\label{Rmatprop}

Trigonometric $U_q(\Dn{n})$-related $\RR$-matrix given by \r{tRalt} possesses following properties.
\begin{itemize}

\item {\it Scaling invariance}
 \begin{equation}\label{scaling}
  \RR(\beta u,\beta v)=\RR(u,v)
 \end{equation}
 for any complex parameter $\beta$ which is not equal to zero.

\item {\it $\xi$-parity}
 \begin{equation}\label{eq: xi parity}
  \left.\RR(u,v) \right|_{\xi \to - \xi} = \RR(u, v),
 \end{equation}
which distinguishes it from $R$-matrix \eqref{D2Rmat}.

\item {\it D-invariance}
 \begin{equation}\label{comm1}
  \Der_1\, \Der_2\, \RR_{12}(u,v)=\RR_{12}(u,v)\, \Der_1\, \Der_2,
 \end{equation}
 where $\Der$ is defined by \r{Der} and $\Der_1=\Der\ot\Id$, $\Der_2=\Id\ot\Der$.

\item {\it U-invariance}
 \begin{equation}\label{comm2}
  \Ut_1\, \Ut_2\, \RR_{12}(u,v)=\RR_{12}(u,v)\, \Ut_1\, \Ut_2\,,
 \end{equation}
where $U$ is given by \r{Ut}.

\item {\it Twist-invariance}

The equation
 \begin{equation}\label{eq:twist}
   \left[ R_{12}(u, v), \theta \otimes \theta \right] = 0
 \end{equation}
 has only two classes of non-degenerated solutions: 1) $\theta = \theta_{d}$, 2) $\theta = \Ut \theta_{d}$,
 where $\theta_{d}$ is diagonal and ${\theta_{d}}^{\mathrm{t}} \theta_{d} \propto \Id$. $D$- and $U$- invarince are particular cases of the twist-invariance.

\item {\it Flip crossing symmetry}
 \begin{equation}\label{reflec1}
  \RR_{12}(-u, v) = \Ut_1\  \RR_{12}(u, v)\  \Ut_1\,.
 \end{equation}

\item {\it Crossing symmetry}
\begin{equation}\label{cross1}
  \RR_{12}(u, v)^{{\rm t}_1}\,  = \TT_1^{-2}\ \Der_1^{-1} \RR_{21}(v \xi, u)\, \Der_1\ \TT_1^2, \quad
  \RR_{12}(u, v)^{{\rm t}_2}\,  = \TT_2^{-2}\ \Der_1^{-1} \RR_{21}(v \xi, u)\, \Der_1\ \TT_2^2\,,
\end{equation}
where $\TT$ is defined by \r{Tm1} and $\TT_1=\TT\ot\Id$, $\TT_2=\Id\ot\TT$.

\item {\it Double transposition symmetry}
 \begin{equation}\label{reflec}
  \RR_{12}(u,v)^{{\rm t}_1{\rm t}_2}=\TT_1^{-2}\ \TT_2^{-2}\ \RR_{12}(u,v)\
  \TT_1^{2}\ \TT_2^{2}\,.
\end{equation}

\item {\it Unitarity}
\begin{equation}\label{unitar}
\RR_{12}(u,v)\cdot \RR_{21}(v,u)=f(u^2,v^2)f(v^2,u^2)\ \Id\ot\Id\,,
\end{equation}
where $\RR_{21}(u,v)=\Per_{12}\,\RR_{12}(u,v)\,\Per_{12}$.

\item {\it Crossing unitarity}
\begin{equation}\label{cross2}
 \Der^2_1\ \RR_{12}(v\xi^2,u)^{{\rm t}_1}\
 \Der_1^{-2}\ \RR_{21}(u,v)^{{\rm t}_1}
 = f(u^2,v^2\xi^2)f(v^2\xi^2,u^2)\ \Id\ot\Id\,.
\end{equation}

\item {\it Yang-Baxter equation}
 \begin{equation}\label{YB}
\RR_{12}(u,v)\cdot \RR_{13}(u,w)\cdot \RR_{23}(v,w)=
\RR_{23}(v,w)\cdot \RR_{13}(u,w)\cdot \RR_{12}(u,v)\,,
\end{equation}
where subscripts of $\RR$-matrices mean the indices of the spaces $\CC^N$ where
they act nontrivially.

\item {\it Pole structure}\\
$\RR$-matrix has simple poles at $u=\pm v$ and $u=\pm v\xi$ with residues
\begin{equation}\label{poles}
\begin{split}
\left.\frac{2 (u-v)}{u (q-q^{-1})}\ \RR_{12}(u,v)\right|_{u=v}&=\Per_{12},\\
\left.\frac{2 (u+v)}{u (q-q^{-1})}\ \RR_{12}(u,v)\right|_{u=-v}&=\Ut_1\ \Per_{12}\ \Ut_1,
\\
\left.\frac{2 (v\xi-u) }{u(q-q^{-1})}\ \RR_{12}(u,v)\right|_{u=v\xi}&=\Der_1^{-1} \ \TT_1^2 \
\Per_{12}^{{\rm t}_1}
\ \TT_1^{-2}\ \Der_1\,, \\
\left.-\frac{2 (v\xi+u) }{u(q-q^{-1})}\ \RR_{12}(u,v)\right|_{u=-v\xi}&=
\Ut_1 \ \Der_1^{-1}\ \TT_1^2\ \Per_{12}^{{\rm t}_1}\ \TT_1^{-2}\ \Der_1\ \Ut_1\,.
\end{split}
\end{equation}
For $n>2$ one can get  from \r{cross1} and \r{cross2}
 \begin{equation}\label{cen2}
 f(\xi^2,1)f(1,\xi^2)\
\left.\frac{2(v\xi^2 -u)}{u(q-q^{-1})}\Big(\RR_{21}(u,v)^{{\rm t}_1}
\Big)^{-1}\right|_{u=v\xi^2}=
 \Der_1^{2}\ \Per_{12}^{{\rm t}_1}\   \Der_1^{-2}\,.
\end{equation}
This pole of $\RR$-matrix will be important to get Nazarov type
central elements \cite{Nazarov1991centers} in the algebra $U_q(\Dn{n})$.

\item {\it Scaling limit}\\
One may check that in the scaling limit $\epsilon\to 0$ and
$u\to e^{\epsilon u}$, $v\to e^{\epsilon v}$, $q\to e^{\epsilon c }$,
trigonometric $\RR$-matrix \r{D2Rmat}   goes into rational $\mathfrak{o}_{2n}$-invariant
$R$-matrix
\begin{equation}\label{R-rat}
\RR(u,v)\ =\  \Id\ot \Id+\frac{c}{u-v}\ \Per\ -
\frac{c }{u-v+c\kappa}\ \Qer\,,
\end{equation}
where
\begin{equation}\label{PQer}
\Qer=\sum_{1\leq i,j\leq N}  \E_{ij} \ot \E_{i'j'}=\Per^{{\rm t}_2}=\Per^{{\rm t}_1}
\end{equation}
and $\kappa=n-1$.

\end{itemize}

The proof of most of the above properties follows directly from the explicit form of the R-matrix \eqref{tRalt}.
Proofs of the crossing relations and the pole structure of $\RR$-matrix \r{tRalt}
can be found in appendix~\ref{sec: Crossing symmetries and pole structure}.

\section{Algebra $U_q(\Dn{n})$ in $\RR$-matrix form}\label{Algebra}

The algebra $U_q(\Dn{n})$ over $\CC(q^{1/2})$ is a collection of the algebraically dependent
generators
$\LL^\pm_{i,j}[\pm m]$, $1\leq i,j\leq 2n$, $m\in \mathbb{Z}_{\ge0}$. These algebraic dependences
are related to the centers in this algebra which will be described in the
section~\ref{sec-cent}.
$U_q(\Dn{n})$ as an associative algebra contains also a unity generator
$\mathbf1$, which satisfy the property $\mathbf1\cdot x=x\cdot\mathbf1=x$.
 Zero modes $\LL^\pm_{i,j}[0]$ plays important role in the algebra
$U_q(\Dn{n})$ and satisfy
\begin{equation}\label{restr}
\LL^+_{j,i}[0]=\LL^-_{i,j}[0]=0, \quad i<j,\quad \LL^+_{i,i}[0]\LL^-_{i,i}[0]=\LL^-_{i,i}[0]\LL^+_{i,i}[0]=
\mathbf1\,.
\end{equation}

As usual, one can collect the generators of the algebra $U_q(\Dn{n})$
into formal generating series using spectral parameter $u$
\begin{equation}\label{series}
\LL^\pm_{i,j}(u)=\sum_{m=0}^\infty \LL^\pm_{i,j}[\pm m] u^{\mp m}
\end{equation}
and combine them into matrices
\begin{equation}\label{Lop}
\LL^\pm(u)=\sum_{i,j=1}^N \E_{ij}\ot \LL^\pm_{i,j}(u)\in {\rm End}(\CC^N)\ot
U_q(\Dn{n})[[u,u^{-1}]]
\end{equation}
which we call {\it fundamental} $\LL$-operators.
The commutation relations in the algebra $U_q(\Dn{n})$ are given
by the  $\RR\LL\LL$ commutation relations in $(\CC^N)^{\ot2}\ot U_q(\Dn{n})[[u,u^{-1}]]$
\begin{equation}\label{RLL}
\RR(u,v) \cdot (\LL^\mu(u) \ot \Id)\cdot (\Id\ot \LL^\rho(v))=
(\Id\ot \LL^\rho(v))\cdot (\LL^\mu(u) \ot \Id)\cdot \RR(u,v)\,,
\end{equation}
where $\mu,\rho=\pm$  and rational functions entering $\RR$-matrix \r{tRalt}
 should be understood as series over $v/u$ for $\mu=+$, $\rho=-$ and
as series over $u/v$ for $\mu=-$, $\rho=+$. For $\mu=\rho$
these rational functions can be either series over the ratio $v/u$ or the ratio  $u/v$.
One  can get
the commutation relations of the algebra $U_q(\Dn{n})$  in terms
of matrix entries \r{series}. We will not present this expression explicitly because
they are rather bulky.
One can also check that restrictions \r{restr} are compatible with the
commutation relations \r{RLL}.

Collections of the generating series $\LL^\mu_{i,j}(u)$ for $\mu=\pm$
form standard Borel subalgebras $U^\pm_q(\tilde{\mathfrak{g}})$ in the quantum
loop algebra $U_q(\Dn{n})$.

\subsection{Centers in $U_q(\Dn{n})$}\label{sec-cent}

To describe centers of the RTT-algebra \eqref{RLL}, we use methods developed in \cite{JingLiuMolev2018Y(BCD)}.
There are two types of the central elements related to the poles of
$\RR$-matrix at $u=-v$ and $u=v\xi$. We denote them as $Z^p_\pm$ and
$Z^c_{\pm}$ respectively. The following proposition gives them explicit expressions.
\begin{prop}\label{centers} The central elements $Z^p_\pm$ and
$Z^c_{\pm}$ are described by the expressions
\begin{equation}\label{Zp}
  Z^p_\pm(u)\ \Id =
  \Ut\  \big( \LL^\pm(u) \big)^{-1}  \Ut\ \LL^\pm(-u)
  =    \LL^\pm(-u)\ \Ut \big( \LL^\pm(u) \big)^{-1} \Ut
\end{equation}
and
\begin{equation}\label{Zc}
    Z_\pm^c(u)\ \Id =
    \Derd\ \LL^\pm(u\xi)^{{\rm t}}\ \Derd^{-1}\ \LL^\pm(u)=
    \LL^\pm(u)\ \Derd\ \LL^\pm(u\xi)^{{\rm t}}\ \Derd^{-1}\,,
\end{equation}
where
\begin{equation}\label{tDer}
{\Derd}=\Der\cdot \TT^{2}\,.
\end{equation}
\end{prop}

\proof To simplify the proof, we skip the superscripts of $\LL$-operators.
Starting from $\RR\LL\LL$ commutation relations \r{RLL} and using
second line in \r{poles} one obtains at  the point $v = - u$
\begin{equation*}
    \Ut_1 \Per_{12} \Ut_1\ \LL_1(u) \LL_2(-u) = \LL_2(-u) \LL_1(u)\ \Ut_1 \Per_{12} \Ut_1\,.
\end{equation*}
Multiplying by $\Ut_1$ from the left and the right and then using the property of the permutation
operator $\Per_{12}$ one gets
\begin{equation*}
     \Ut_2\ \LL_2(u)\ \Ut_2\ \LL_1(-u) = \LL_2(-u)\ \Ut_1\ \LL_1(u)\ \Ut_1 \,.
\end{equation*}
After conjugation, we have
\begin{equation*}
   \left( \LL_2(-u)\ \right)^{-1}  \Ut_2\ \LL_2(u)\ \Ut_2 =  \Ut_1\ \LL_1(u)\ \Ut_1  \left( \LL_1(-u) \right)^{-1} \,.
\end{equation*}
The latter equation proves that operators
\begin{equation*}
  \Ut\  \big( \LL(u) \big)^{-1}  \Ut\ \LL(-u)  =    \LL(-u)\ \Ut \big( \LL(u) \big)^{-1} \Ut=
   Z^p(u)\ \Id
\end{equation*}
are proportional to the unity operator in ${\rm End}(\CC^N)$.

Let us prove that $Z^p(u)$ is a center in the algebra $U_q(\Dn{n})$
\begin{equation*}
\begin{split}
    Z^p(u) \LL_2(v)&=Z^p(u)\ \Id_1 \LL_2(v) =
    \Ut_1\  \big( \LL_1(u) \big)^{-1}  \Ut_1\ \LL_1(-u)  \LL_2(v) = \\
    &=\Ut_1\  \big( \LL_1(u) \big)^{-1}  \Ut_1\ \big(\RR_{12}(-u,v)\big)^{-1}
    \LL_2(v) \LL_1(-u) \RR_{12}(-u,v) = \\
    &=\Ut_1\  \big( \LL_1(u) \big)^{-1}  \big(\RR_{12}(u,v)\big)^{-1}
    \LL_2(v)\ \Ut_1\  \LL_1(-u) \RR_{12}(-u,v) = \\
    &=\LL_2(v)\ \Ut_1\  \big(\RR_{12}(u,v)\big)^{-1} \big( \LL_1(u) \big)^{-1} \
    \Ut_1\  \LL_1(-u) \RR_{12}(-u,v) = \\
    &=\LL_2(v)  \big(\RR_{12}(-u,v)\big)^{-1} \Ut_1\
    \big( \LL_1(u) \big)^{-1} \ \Ut_1\  \LL_1(-u) \RR_{12}(-u,v) = \\
    &=\LL_2(v)  \big(\RR_{12}(-u,v)\big)^{-1} Z^p(u) \, \Id_{1} \ \RR_{12}(-u,v) =
    \LL_2(v) Z^p(u).
    \end{split}
\end{equation*}

Considering $\RR\LL\LL$ relation
at the point $u = \xi v$ one can obtain using third line in \r{poles}
\begin{equation*}
    \Der_1^{-1}\ \TT_1^2\ \Per^{\rm t_{1}}_{12}\ \TT_1^{-2}\ \Der_1\ \LL_1(v \xi) \LL_2(v) =
    \LL_2(v) \LL_1(v\xi)\ \Der_1^{-1}\ \TT_1^2\ \Per^{\rm t_{1}}_{12}\ \TT_1^{-2}\ \Der_1
\end{equation*}
and then prove that the operators
\begin{equation}\label{cent-eq}
    \Derd\ \LL(v\xi)^{{\rm t}}\ \Derd^{-1}\ \LL(v)=
    \LL(v)\ \Derd\ \LL(v\xi)^{{\rm t}}\ \Derd^{-1}=  Z^c(v)\ \Id
\end{equation}
are proportional to the unity operator in ${\rm End}(\CC^N)$.
One should take into account that $\Der^{\rm t}=\Der^{-1}$.
The centrality of  $Z^c(u)$ can be proved using \r{RLL} as above. \qed

\begin{remark}\label{rm: zpzp}
  The following relation is valid
\begin{equation}\label{eq: zpzp}
 Z_{\pm}^{p} (u) Z_{\pm}^{p}(-u) = 1,
\end{equation}
which is a direct consequence of the equation \eqref{Zp}.
\end{remark}

We use \r{Zp} and \r{Zc}
to find the set of algebraically independent generators in terms of the
Gaussian coordinates of the $\LL$-operators.

Besides central elements given by Proposition~\ref{centers}, these are
Nazarov type central elements $Z^{\rm N}_\pm(v)$
\cite{Nazarov1991centers} due to the
equality \r{cen2}.
\begin{prop}\label{Nazcen}
There are central elements $Z^{\rm N}_\pm(v)$ in the algebra $U_q(\Dn{n})$
\begin{equation}\label{ZN1}
Z^{\rm N}_\pm(v)\ \Id=
\Der^2\ \LL^\pm(v\xi^2)^{\rm t}\ \Der^{-2}\ \sk{\LL^\pm(v)^{-1}}^{\rm t}=
\sk{\LL^\pm(v)^{-1}}^{\rm t}\ \Der^2\ \LL^\pm(v\xi^2)^{\rm t}\ \Der^{-2}\,.
\end{equation}
\end{prop}
The proof of this proposition is the same as in the paper
\cite{LiashykPakuliak2022U_q(ABCD)}.

\begin{prop}
Assuming that the central elements \r{Zc} are invertible they
are related to the central elements \r{ZN1} by the relation
\begin{equation}\label{ZN2}
Z^{\rm N}_\pm(v)=Z^{c}_\pm(v\xi)Z^{c}_\pm(v)^{-1}\,.
\end{equation}
\end{prop}

For the proof of this proposition, it is sufficient to consider
\begin{equation*}
Z_\pm^c(v\xi)\ \Id =
    \Der\ \TT^2\ \LL^\pm(v\xi^2)^{{\rm t}}\ \TT^{-2}\ \Derd^{-1}\ \LL^\pm(v\xi)
\end{equation*}
and express from \r{Zc} $\LL$-operator
\begin{equation*}
\LL^\pm(v\xi)=Z^c_\pm(v)\ \Der^{-1}\ \TT^2\ \Big(\LL^\pm(v)^{-1}\Big)^{\rm t}\TT^{-2}\ \Der\,.
\end{equation*}
Substituting the second equality into the first one, we observe that twisting matrix $\TT$
is canceled, and one gets
\begin{equation*}
\begin{split}
Z_\pm^c(v\xi)\ \Id&=Z^c_\pm(v)\ \Der^{-1}\ \TT^2\ \Big(
\Der^2\ \LL^\pm(v\xi^2)^{\rm t}\ \Der^{-2}\ \sk{\LL^\pm(v)^{-1}}^{\rm t}
\Big)\ \TT^{-2}\ \Der=\\
 &=Z^c_\pm(v)\ Z^{\rm N}_\pm(v)\, \Id.
 \end{split}
\end{equation*}
\qed

\subsection{Gaussian coordinates}\label{GCsect}

One can introduce Gaussian coordinates of the fundamental
$\LL$-operators in several different ways. There are two preferable choices of
the Gaussian decomposition, which are suitable for the further application of the
Gaussian coordinates in the algebraic Bethe ansatz. Either one decomposes
the fundamental $\LL$-operators in the product of the upper triangular, diagonal and
lower triangular matrices or one does the same for the transposed $\LL$-operators
$\LL^\pm(u)^{\rm t}$. Different Gaussian decompositions allow us to obtain the
different recurrence relations for the off-shell Bethe vectors (see for example
the papers \cite{OskinPakuliakSilantyev2009recrel} and
\cite{HutsalyukLiashykPakuliakRagoucySlavnov2017current}, where different Gaussian
decompositions were used to get such recurrence relations for the off-shell Bethe
vectors in the quantum integrable models associated with the quantum affine algebra
$U_q(\mathfrak{gl}_N)$ and the super-Yangian double $DY(\mathfrak{gl}(m|n))$
respectively.)

In this paper, we will explore the Gaussian decomposition for the transposed
$\LL$-operators  $\LL^\pm(u)^{\rm t}$  given by
the products of upper triangular $\mathbf{F}^\pm(u)$, diagonal
$\mathbf{K}^\pm(u)$ and lower triangular  $\mathbf{E}^\pm(u)$ matrices
\begin{equation}\label{Gauss3}
\begin{split}
&\LL^\pm(u)^{\rm t}= \sum_{1\leq i,j\leq 2n} \LL^\pm_{i,j}(u)\ \E^{\rm t}_{ij}=
\mathbf{F}^\pm(u)^{{\rm t}}\times \mathbf{K}^\pm(u)^{{\rm t}}\times \mathbf{E}^\pm(u)^{{\rm t}}=\\
&\quad=\sk{\sum_{1\leq i\leq j\leq 2n} \E_{ij}^{\rm t}\ \FF^\pm_{j,i}(u)}
\times \sk{\sum_{1\leq i\leq 2n} \E_{ii}^{\rm t}\ k^\pm_{i}(u)}\times
\sk{\sum_{1\leq i\leq j\leq 2n} \E_{ji}^{\rm t}\ \EE^\pm_{i,j}(u)}\,,
\end{split}
\end{equation}
where the Gaussian coordinates $\FF^\pm_{i,i}(u)=\EE^\pm_{i,i}=\mathbf{1}$ and
 $\FF^\pm_{n+1,n}(u)$ and $\EE^\pm_{n,n+1}(u)$ vanish.
This vanishing is similar to the case of the quantum loop algebra $U_q(D^{(1)}_n)$
noted for the first time in \cite{JingLiuMolev2020Uq(BD)}.
One can prove this vanishing by analyzing the $\RR\LL\LL$-relation, but in the next
subsection, we will provide arguments explaining this fact using notion of normal
ordering of the Gaussian coordinates.

 Using multiplication rule
 \begin{equation*}\label{muru}
\E^{\rm t}_{ij}\cdot \E^{\rm t}_{kl}=\delta_{il}\ \E^{\rm t}_{kj}
\end{equation*}
one can obtain explicit form of the matrix entries $\LL^\pm_{i,j}(u)$ in terms
of the introduced Gaussian coordinates
\begin{equation}\label{Gauss1}
\LL^\pm_{i,j}(u)=\sum_{1\leq\ell\leq{\rm min}(i,j)}
\FF^\pm_{j,\ell}(u)\ k^\pm_{\ell}(u)\ \EE^\pm_{\ell,i}(u)\,.
\end{equation}

Matrix entries of $\LL$-operators and Gaussian coordinates are not algebraically independent
due to the existance of  constraints \eqref{Zp} and \eqref{Zc}.

Consider the solution of the equation
\begin{equation}\label{eq: Zp ff}
  Z^{p}_{\pm}(u) = \frac{f^{\pm}(u)}{f^{\pm}(-u)},
\end{equation}
which agrees with \eqref{eq: zpzp}.
Using the freedom to rescale $L^{\pm}$-operators in the RLL relation \eqref{RLL}
\begin{equation}\label{eq: L rescale}
  L^{\pm} \to f^{\pm}(u) L^{\pm}(u)
\end{equation}
we can fix the value of
\begin{equation}\label{eq: zp 1}
Z^p_\pm(u)=\mathbf{1}.
\end{equation}
We also denote
\begin{equation}\label{eq: zc z}
Z^c_\pm(u)=z^\pm(\xi u),
\end{equation}
where $z^\pm(u)$ are series with respect to non-negative
powers of $u^{-1}$ and $u$ respectively and coefficients belong to the center of RLL algebra \eqref{RLL}.

Fixation of $Z^p_\pm(u)$ to the unity \eqref{eq: zp 1} is equivalent
to the following relations
\begin{equation}\label{D22-1}
\LL^\pm(-u)=\Ut\ \LL^\pm(u)\ \Ut\,.
\end{equation}
As a direct consequence of the relation \r{D22-1} and Gaussian decomposition
\r{Gauss3} one gets
\begin{equation}\label{Dpm}
\begin{split}
\FF_{j,n+1}(u)=\FF_{j,n}(-u),\quad &\EE_{n+1,j}(u)=\EE_{n,j}(-u),\quad j>n+1\,,\\
\FF_{n+1,i}(u)=\FF_{n,i}(-u),\quad &\EE_{i,n+1}(u)=\EE_{i,n}(-u),\quad i<n\,,\\
k_{n+1}(u)&=k_n(-u)\,.
\end{split}
\end{equation}
On the other hand, according to \r{D22-1}
the Gaussian coordinates are the even series of the spectral parameters
$\FF^\pm_{j,i}(u)=\FF^\pm_{j,i}(-u)$, $\EE^\pm_{i,j}(u)=\EE^\pm_{i,j}(-u)$
and $k^\pm_j(u)=k^\pm_j(-u)$ for $i\not=n,n+1$ and $j\not=n,n+1$.

On the other hand, using notation for the second central element $Z^c_\pm(u)$ \eqref{eq: zc z}
leads to the following relations between $\LL$-operators
\begin{equation*}\label{ident}
\Derd\ \LL^\pm(\xi u)^{{\rm t}}\ \Derd^{-1}= z^\pm(\xi u) \LL^\pm(u)^{-1}\,,
\end{equation*}
or
\begin{equation}\label{identalt}
\Derd\ \hLL^\pm(\xi u)\ \Derd^{-1}= z^\pm(\xi u)^{-1} \LL^\pm(u)\,,
\end{equation}
where transposed-inversed $\LL$-operators $\hLL^\pm(u)$ are defined as
\begin{equation}\label{hLL}
\hLL^\pm(u)=\Big(\LL^\pm(u)^{{\rm t}}\Big)^{-1}\,.
\end{equation}

Note that due to the commutativity $[\Derd,\Ut]=0$ and relation \eqref{D22-1} the central element
$z^\pm(u)=z^\pm(-u)$ is an even series of the formal parameter $u$ or $u^{-1}$.

Gaussian decomposition \r{Gauss3} allows to obtain similar decomposition
of $\LL$-operators $\hLL^\pm(u)$. Indeed,
 taking the inverse of both sides of the equality  \r{Gauss3}
\begin{equation}\label{Gau33}
\begin{split}
\hLL^\pm(u)&= \Big(\big(\LL^\pm(u)\big)^{\rm t}\Big)^{-1}
= \sum_{1\leq i,j\leq 2n} \E_{ij}\ \hLL^\pm_{i,j}(u)=\\
&=\sk{\sum_{1\leq i\leq j\leq 2n} \E_{ji}^{\rm t}\ \tEE^\pm_{i,j}(u)}\times
\sk{\sum_{1\leq i\leq 2n} \E_{ii}^{\rm t}\ k^\pm_{i}(u)^{-1}}\times
\sk{\sum_{1\leq i\leq j\leq 2n} \E_{ij}^{\rm t}\ \tFF^\pm_{j,i}(u)}
\end{split}
\end{equation}
one gets
 \begin{equation}\label{GChL}
 \hLL^\pm_{i,j}(u)=
 \sum_{ \ell\leq {\rm min}(i,j)}\tEE^\pm_{i',\ell'}(u)\ k^\pm_{\ell'}(u)^{-1}\
 \tFF^\pm_{\ell',j'}(u)\,.
 \end{equation}

Gaussian coordinates $\tFF^\pm_{j,i}(u)$ and  $\tEE^\pm_{i,j}(u)$
in \r{Gau33} and \r{GChL} satisfy the relations
\begin{equation}\label{invrel}
\sum_{i\leq\ell\leq j}\FF^\pm_{j,\ell}(u)\tFF^\pm_{\ell,i}(u)= \delta_{ij},\quad
\sum_{i\leq\ell\leq j}\EE^\pm_{i,\ell}(u)\tEE^\pm_{\ell,j}(u)=\delta_{ij}\,.
\end{equation}
These equations for the  Gaussian coordinates $\tFF^\pm_{j,i}(u)$ and $\tEE^\pm_{i,j}(u)$
can be resolved to get
\begin{align}
\tFF^\pm_{j,i}(u) &=  -\FF^\pm_{j,i}(u)+\sum_{\ell=1}^{j-i-1}(-)^{\ell+1}
\sum_{j>i_\ell>\cdots>i_1>i} \FF^\pm_{j,i_\ell}(u)  \FF^\pm_{i_\ell,i_{\ell-1}}(u)\cdots
\FF^\pm_{i_2,i_1}(u)  \FF^\pm_{ i_1,i}(u)\,, \label{tFF}  \\
\tEE^\pm_{i,j}(u) &=  -\EE^\pm_{i,j}(u)+\sum_{\ell=1}^{j-i-1}(-)^{\ell+1}
\sum_{j>i_\ell>\cdots>i_1>i} \EE^\pm_{i,i_1 }(u)  \EE^\pm_{i_1,i_2}(u)\cdots
\EE^\pm_{i_{\ell-1},i_\ell}(u)  \EE^\pm_{i_\ell,j}(u)\,. \nonumber 
\end{align}

\subsection{Normal ordering of the Gaussian coordinates}\label{nor-ord}

It is known \cite{KHOROSHKIN1993445} that the modes of Gaussian
coordinates coincide with Cartan-Weyl generators in the quantum affine algebras.
The set of Cartan-Weyl generators has natural {\it convex} ordering
\cite{Beck1994ConvexBases}, which was used in
\cite{EnriquezKhoroshkinPakuliak2007Projections} to define
the normal ordering of the Gaussian coordinates. Let us shortly
remind this construction.

Let $U^\pm_f$, $U^\pm_e$ and $U^\pm_k$ be subalgebras of $U_q(\gaf)$
formed by the modes of the Gaussian coordinates $\FF^\pm_{j,i}(u)$, $\EE^\pm_{i,j}(u)$
and $k^\pm_j(u)$.  Convex ordering of the Cartan-Weyl generators
\cite{Beck1994ConvexBases} implies  the normal ordering of the subalgebras
formed by the Gaussian coordinates. These orderings are
\begin{equation}\label{order1}
\cdots\prec U^-_k\prec U^-_f\prec U^+_f\prec U_k^+\prec U^+_e\prec U^-_e\prec U^-_k\prec\cdots
\end{equation}
or
\begin{equation}\label{order2}
\cdots\prec U^+_k\prec U^+_f\prec U^-_f\prec U_k^-\prec U^-_e\prec U^+_e\prec U^+_k\prec\cdots
\end{equation}
This ordering was called convex because
if one places subalgebras $U^\pm_f$, $U^\pm_e$ and $U^\pm_k$
onto circles
\begin{equation*}
\begin{array}{ccccc}
&U^-_e&&U^+_e& \\[3mm]
U^-_k& &\circlearrowleft& & U^+_k\\[3mm]
&U^-_f&&U^+_f&
\end{array}
\qquad\mbox{or} \qquad
\begin{array}{ccccc}
&U^-_e&&U^+_e& \\[3mm]
U^-_k& &\circlearrowright& & U^+_k\\[3mm]
&U^-_f&&U^+_f&
\end{array}
\end{equation*}
then one can prove  that in both types
of ordering \r{order1} and  \r{order2}
the unions of subalgebras $U^\pm_f$, $U^\pm_e$ and $U^\pm_k$
along the smallest arcs between starting and ending points are
subalgebras in the quantum loop algebra  $U_q(\gaf)$.
For example, the union of subalgebras $U^+_f\cup U^+_k$
or $U^-_f\cup U^+_f\cup U^+_k$ or  $U_q^+(\gaf)=U^+_f\cup U^+_k\cup U^+_e$
and so on are subalgebras in  $U_q(\gaf)$.

This normal ordering yields a powerful practical tool to
get relations for the Gaussian coordinates. In any relation
which contains Gaussian coordinates of the different types one has to
order first all monomials according to \r{order1} or \r{order2}  and then single out
all the terms that belong to the one of subalgebras which is composed of the
Gaussian coordinates of the necessary type.
Using this procedure, one can get quite easy the commutation relations
between Gaussian coordinates from $\RR\LL\LL$-commutation
relations \r{RLL} \cite{LiashykPakuliak2022U_q(ABCD)}.
We call this procedure {\it a restriction} to subalgebras in $U_q(\gaf)$.

Subalgebras $U^\pm_q(\gaf)=U^\pm_f\cup U^\pm_k\cup U^\pm_e$ were already introduced above as  Borel subalgebras in $U_q(\gaf)$.
To describe the so-called 'new realization' of these algebras in terms of the currents
\cite{Drinfeld1988NewRealisation}
one has to consider different type Borel subalgebras $U_f=U^-_f\cup U^+_f\cup U^+_k$
and $U_e=U^+_e\cup U^-_e\cup U^-_k$.

Considering intersections between different types of Borel subalgebras
one can introduce    projections $\Pfpm$  and $\Pepm$
which  were  investigated in
\cite{EnriquezKhoroshkinPakuliak2007Projections} for the ordering
\r{order1}.  One can check that
 the action of the projections $\Pfpm$  and $\Pepm$
onto Borel subalgebras $U_f$ and $U_e$ described in this paper
coincides with restrictions on subalgebras
$U_f^\pm$ and $U_e^\pm$  defined for the ordering \r{order1}.

Assume now that $\FF^\pm_{n+1,n}(u)\not=0$ and $\EE^\pm_{n,n+1}(u)\not=0$.
Then according to \r{Gauss1} matrix entries $\LL^+_{n,n+1}(u)$ and
$\LL^+_{n+1,n}(u)$ take the form
\begin{equation*}
\begin{split}
\LL^+_{n,n+1}(u)&=\FF^+_{n+1,n}(u)k^+_n(u)+ \sum_{1\leq\ell\leq n-1}
\FF^+_{n+1,\ell}(u)\ k^+_{\ell}(u)\ \EE^+_{\ell,n}(u)\,,\\
\LL^+_{n+1,n}(u)&=k^+_n(u)\EE^\pm_{n,n+1}(u)+ \sum_{1\leq\ell\leq n-1}
\FF^+_{n,\ell}(u)\ k^+_{\ell}(u)\ \EE^+_{\ell,n+1}(u)\,.
\end{split}
\end{equation*}
Note that Gaussian coordinates in these expressions are normal ordered according
to the ordering \r{order1}.
On the other hand, these matrix elements are related due to \r{D22-1}
as follows
\begin{equation}\label{van-eq}
\LL^+_{n,n+1}(u)=\LL^+_{n+1,n}(-u)\,.
\end{equation}
Considering restriction of the latter equality to the subalgebra $U_f^+\cup U^+_k$
one obtains that $\FF^+_{n+1,n}(u) \cdot k^+_n(u)=0$ which results in vanishing
of the Gaussian coordinates $\FF^+_{n+1,n}(u)$ due to invertibility of the current
$k^+_n(u)$. Analogously, considering restriction of \r{van-eq} to the subalgebra
$ U^+_k\cup U_e^+$ one gets that $\EE^+_{n,n+1}(u)=0$. Analogous
to \r{van-eq} equality for the 'negative'
matrix entries results to the same vanishing of $\FF^-_{n+1,n}(u)$
and $\EE^-_{n,n+1}(u)$ since these matrix entries are also normal ordered
with respect to another ordering \r{order2} and one can restrict them
to the corresponding subalgebras to prove this vanishing.

These arguments will be justified in Section~\ref{cur-real} by the explicit calculation
of the monodromy matrix elements in terms of the Gaussian coordinates.

\section{Embedding $U_q(\Dn{n-1})\hookrightarrow  U_q(\Dn{n})$ }\label{embed}

In this sections we will investigate embedding
$U_q(\Dn{n-1})\hookrightarrow  U_q(\Dn{n})$ in the way similar to the paper \cite{JingLiuMolev2018Y(BCD)}.
We will write $\RR^{(n)}(u,v)$, $\Per^{(n)}(u,v)$, $\Qer^{(n)}(u,v)$, etc.,
to associate corresponding matrices with the algebra $U_q(\Dn{n})$ and
$\RR^{(n-1)}(u,v)$, $\Per^{(n-1)}(u,v)$, $\Qer^{(n-1)}(u,v)$, etc., with
the algebra $U_q(\Dn{n-1})$.
Let us define $\LL$-operator $\MM(u)$ with matrix entries
$\MM^\pm_{i,j}(u)$ for  $1<i,j<2n$ and Gaussian decomposition
\begin{equation}\label{Gaussian2}
\MM^\pm_{i,j}(u)=\sum_{2\leq\ell\leq{\rm min}(i,j)}
\FF^\pm_{j,\ell}(u)\ k^\pm_{\ell}(u)\ \EE^\pm_{\ell,i}(u)=
\LL^\pm_{i,j}(u)-\LL^\pm_{1,j}(u)\LL^\pm_{1,1}(u)^{-1}\LL^\pm_{i,1}(u).
\end{equation}

Let $|\ell\>$ for $\ell=1,\ldots,2n$ be a basis in $\CC^{2n}$. We denote by $|i,j\>=
|i\>\ot |j\>$ the basis elements in $\CC^{2n}\ot\CC^{2n}$. Analogously
we define $|i,j,\ell\>\in \big(\CC^{2n}\big)^{\ot3}$ and
$|i,j,\ell,m\>\in \big(\CC^{2n}\big)^{\ot4}$. Dual vectors
$\<\ell|$, $\<i,j|$, $\<i,j,\ell|$ and  $\<i,j,\ell,m|$ are defined similarly.

Let $\Lf^{(1,2)}(u)$ be a fused $\LL$-operator\footnote{To avoid
bulky notations we skip superscripts of $\LL$-operators in this section.}
\begin{equation}\label{Lf}
\Lf^{(1,2)}(u)=\RR_{12}(1,q)\LL^{(1)}(u)\LL^{(2)}(qu)=\LL^{(2)}(qu)\LL^{(1)}(u)\RR_{12}(1,q)\,.
\end{equation}
Let us calculate $(i,j;1,1)$ matrix element of this fused $\LL$-operator.
One has
\begin{multline}\label{e2}
  \Lf_{i,j;1,1}(u)=\<i,1|\Lf^{(1,2)}(u)|j,1\>=\\
  \begin{cases}
\LL_{i,j}(u)\LL_{1,1}(qu)-q\LL_{1,j}(u)\LL_{i,1}(qu)\quad&\mbox{for}\quad i\not=n,n+1\,,\\[1mm]
\LL_{1,1}(qu)\LL_{i,j}(u)-q^{-1}\LL_{1,j}(qu)\LL_{i,1}(u)\quad&\mbox{for}\quad j\not=n,n+1\,,\\[1mm]
\LL_{i,j}(u)\LL_{1,1}(qu)-\frac{q+1}{2}\LL_{1,j}(u)\LL_{i,1}(qu)
-\frac{q-1}{2}\LL_{1,j}(u)\LL_{i',1}(qu)
\quad&\mbox{for}\quad i=n,n+1\,,\\[1mm]
\LL_{1,1}(qu)\LL_{i,j}(u)-\frac{q^{-1}+1}{2}\LL_{1,j}(qu)\LL_{i,1}(u)
-\frac{q^{-1}-1}{2}\LL_{1,j'}(qu)\LL_{i,1}(u)
\quad&\mbox{for}\quad j=n,n+1\,,
\end{cases}
\end{multline}
where equation \r{em3} and \r{emm4} have been used.
Using commutation relations
\begin{equation}\label{e3}
\begin{split}
\LL_{i,1}(u)\LL_{1,1}(qu)&=q\LL_{1,1}(u)\LL_{i,1}(qu),\quad i\not=n,n+1\,,\\
\LL_{1,1}(qu)\LL_{1,j}(u)&=q^{-1}\LL_{1,j}(qu)\LL_{1,1}(u),\quad j\not=n,n+1\,,\\
\LL_{i,1}(u)\LL_{1,1}(qu)&=\LL_{1,1}(u)
\sk{\frac{q+1}{2}\LL_{i,1}(qu)
+\frac{q-1}{2}\LL_{i',1}(qu)},\quad i=n,n+1\,,\\
\LL_{1,1}(qu)\LL_{1,j}(u)&=
\sk{\frac{q^{-1}+1}{2}\LL_{1,j}(qu)
+\frac{q^{-1}-1}{2}\LL_{1,j'}(qu)}\LL_{1,1}(u),\quad j=n,n+1
\end{split}
\end{equation}
one can obtain
\begin{equation}\label{e4}
\MM_{i,j}(u)=\Lf_{i,j;1,1}(u)\ \LL_{1,1}(qu)^{-1}=\LL_{1,1}(qu)^{-1}\ \Lf_{i,j;1,1}(u)
 \quad \forall\ 1<i,j<2n\,.
\end{equation}

\begin{lemma}\label{5.2}
There is a commutativity relation
\begin{equation}\label{e5}
\LL_{1,1}(u)\ \MM_{i,j}(v)= \MM_{i,j}(v)\ \LL_{1,1}(u) \,.
\end{equation}
\end{lemma}

To prove this lemma, one has to consider the commutation relations
between $\LL$-operators $\LL(u)$ and $\Lf(v)$
\begin{equation}\label{e6}
\RR_{12}(u,v)\RR_{13}(u,qv)\ \LL^{(1)}(u)\ \Lf^{(2,3)}(v)=
\Lf^{(2,3)}(v)\ \LL^{(1)}(u)\ \RR_{13}(u,qv)\RR_{12}(u,v)
\end{equation}
and take the matrix elements of this commutation relation between vectors
$\<1,i,1|$ and $|1,j,1\>$ for $1<i,j<2n$. Then using equalities \r{em5} proved
in lemma~\ref{5.2l} one obtains
\begin{equation}\label{e7}
\<1,i,1| \LL^{(1)}(u)\ \Lf^{(2,3)}(v) |1,j,1\> = \<1,i,1|  \Lf^{(2,3)}(v)\ \LL^{(1)}(u) |1,j,1\>,
\end{equation}
which implies the statement of lemma~\ref{5.2}. \qed

\subsection{Embedding relations for $\RR(u,v)$}

We will need following formulas for the action of the operator
$\RR(1,q)$ onto specific vectors in $\CC^{2n}\ot\CC^{2n}$.

\begin{equation}\label{em1}
\begin{array}{l}
\RR(1,q)|\ell,\ell\>=0\\[3mm]
\<\ell,\ell|\RR(1,q)=0
\end{array}
\quad\mbox{for}\quad 1\leq \ell<n\quad\mbox{and}\quad
n+1<\ell\leq 2n,
\end{equation}
and
\begin{equation}\label{em3}
\begin{split}
\RR(1,q)|1,\ell\>&=|1,\ell\>-q|\ell,1\>,\quad \RR(1,q)|\ell,1\>=|\ell,1\>-q^{-1}|1,\ell\>\,,\\
\<\ell,1|\RR(1,q)&=\<\ell,1|-q\<1,\ell|,\quad \<1,\ell|\RR(1,q)=\<1,\ell|-q^{-1}\<\ell,1|
\end{split}
\end{equation}
for $1<\ell<n$ and $n+1<\ell<2n$. For $\ell=n,n+1$ the right and left actions of $\RR(1,q)$
to the vectors $|1,\ell\>$, $|\ell,1\>$ and $\<1,\ell|$, $\<\ell,1|$ are
\begin{equation}\label{emm4}
\begin{split}
\RR(1,q)\ |1,\ell\>&=|1,\ell\>-\frac{q+1}{2}|\ell,1\>-\frac{q-1}{2}|\ell',1\>,\\
\RR(1,q)\ |\ell,1\>&=|\ell,1\>-\frac{q^{-1}+1}{2}|1,\ell\>-\frac{q^{-1}-1}{2}|1,\ell'\>\,,\\
\<\ell,1|\ \RR(1,q)&=\<\ell,1|-\frac{q+1}{2}\<1,\ell|-\frac{q-1}{2}\<1,\ell'|\,,\\
\<1,\ell|\ \RR(1,q)&=\<1,\ell|-\frac{q^{-1}+1}{2}\<\ell,1|-\frac{q^{-1}-1}{2}\<\ell',1|\,.
\end{split}
\end{equation}
It follows from \r{emm4} that
\begin{subequations}\label{em4}
\begin{equation}\label{em4a}
\RR(1,q)\Big(|1,\ell\>+|1,\ell'\>\Big)=\Big(|1,\ell\>+|1,\ell'\>\Big)-q
\Big(|\ell,1\>+|\ell',1\>\Big),
\end{equation}
\begin{equation}\label{em4b}
\RR(1,q)\Big(|\ell,1\>+|\ell',1\>\Big)=\Big(|\ell,1\>+|\ell',1\>\Big)-q^{-1}
\Big(|1,\ell\>+|1,\ell'\>\Big),
\end{equation}
\begin{equation}\label{em4c}
\RR(1,q)\Big(|1,\ell\>-|1,\ell'\>\Big)=\Big(|1,\ell\>-|1,\ell'\>\Big)-
\Big(|\ell,1\>-|\ell',1\>\Big),
\end{equation}
\begin{equation}\label{em4d}
\RR(1,q)\Big(|\ell,1\>-|\ell',1\>\Big)=\Big(|\ell,1\>-|\ell',1\>\Big)-
\Big(|1,\ell\>-|1,\ell'\>\Big)\,,
\end{equation}
\end{subequations}
and similar formulas for the left action of $\RR(1,q)$ to the vectors
$\<1,\ell|\pm\<1,\ell'|$ and $\<\ell,1|\pm\<\ell',1|$.
Formulas \r{em1}--\r{em4} follow from the explicit expression \r{RJmat}
and specializations
\begin{equation*}
f(u^2,v^2)\Big|_{v=qu}=0,\quad \gle(u^2,v^2)\Big|_{v=qu}=-q^{-1},\quad
\gri(u^2,v^2)\Big|_{v=qu}=-q,
\end{equation*}
\begin{equation*}
\gle(\pm u,v)\Big|_{v=qu}=-(q^{-1}\pm 1),\quad
\gri(\pm u,v)\Big|_{v=qu}=-(q\pm 1).
\end{equation*}

Using these formulas one can prove the following lemma.
\begin{lemma}\label{5.2l}
\begin{subequations}\label{em5}
There are equalities
\begin{equation}\label{em5a}
\RR_{23}(1,q)\RR_{13}(u,qv)\RR_{12}(u,v)|1,\ell,1\>=f(u^2,v^2)\RR_{23}(1,q)|1,\ell,1\>
\end{equation}
\begin{equation}\label{em5b}
\<1,\ell,1|\RR_{12}(u,v) \RR_{13}(u,qv) \RR_{23}(1,q)=f(u^2,v^2)\<1,\ell,1|\RR_{23}(1,q)
\end{equation}
\end{subequations}
for  $1<\ell<2n$.
\end{lemma}

Proof of both equalities \r{em5} is similar.
We prove only \r{em5a} and  consider separately two cases
$\ell\not= n,n+1$ and $\ell=n,n+1$. We start from the first case.

A direct consequence of \r{em3} is that
\begin{equation}\label{em6}
\RR(1,q)|1,\ell\>=-q\RR(1,q)|\ell,1\>\quad\mbox{for}\quad 1<\ell<2n\quad\mbox{and}
\quad \ell\not=n,n+1.
\end{equation}

We will use presentation of $\RR$-matrix $\RR(u,v)$ in the form \r{tRalt}
written in the form
\begin{equation}\label{RPQform}
\RR(u,v)=\Id\ot\Id+\PPd(u,v)+\QQd(u,v)\,,
\end{equation}
where $\PPd(u,v)$ is given by \r{PP14} and
\begin{equation}\label{QQ14}
\QQd(u,v)=\frac12\Big(\Qer(u,v)+(\Ut\ot\Id)\cdot\Qer(-u,v)\cdot(\Ut\ot\Id)\Big)
\end{equation}
with $\Qer(u,v)$ given by \r{Qeven}.

 According to these formulas
only unity and matrix $\PPd(u,v)$ will contribute to the action
$\RR_{12}(u,v)$ to the vector $|1,\ell,1\>$ and the action
$\RR_{13}(u,qv)$ to the vectors $|1,\ell,1\>$ and $|\ell,1,1\>$. This implies that
\begin{subequations}\label{em7}
\begin{equation}\label{em7a}
\RR_{12}(u,v)|1,\ell,1\>=|1,\ell,1\>+\frac12\Big(\pf_{\ell1}(u,v)+ \pf_{\ell1}(-u,v) \Big)
|\ell,1,1\>
\end{equation}
and
\begin{equation}\label{em7b}
\RR_{13}(u,qv)|1,\ell,1\>=|1,\ell,1\>+\frac12\Big(\pf_{11}(u,qv)+ \pf_{11}(-u,qv) \Big)
|1,\ell,1\>,
\end{equation}
\begin{equation}\label{em7c}
\RR_{13}(u,qv)|\ell,1,1\>=|\ell,1,1\>+\frac12\Big(\pf_{1\ell}(u,qv)+ \pf_{1\ell}(-u,qv) \Big)
|1,1,\ell\>.
\end{equation}
\end{subequations}
Calculating subsequently the action of $\RR$-matrices to the vector $|1,\ell,1\>$
and using \r{em1} and \r{em6} one gets that left hand side of equality \r{em5}
is equal to
\begin{equation*}
\begin{split}
&\Big(1+\frac12\Big(\pf_{11}(u,qv)+ \pf_{11}(-u,qv)\Big)-\\
&\quad-\frac{q}{4}\Big(\pf_{\ell1}(u,v)+ \pf_{\ell1}(-u,v) \Big)
\Big(\pf_{1\ell}(u,qv)+ \pf_{1\ell}(-u,qv)\Big)\Big)\RR_{23}(1,q)|1,\ell,1\>=\\
&=\Big(f(u^2,q^2v^2)-q\gri(u^2,v^2)\gle(u^2,q^2v^2)\Big)\RR_{23}(1,q)|1,\ell,1\>
=f(u^2,v^2)\RR_{23}(1,q)|1,\ell,1\>
\end{split}
\end{equation*}
which proves \r{em5a}. In this calculation we used definitions of the functions $\pf_{ij}(u,v)$
\r{p-fun} and equalities \r{sum-eq}.

Proof of equality \r{em5} for $\ell=n,n+1$ is a little bit more tricky but straightforward.
Using again  \r{tRalt}, one can get that
\begin{equation}\label{em8}
\begin{split}
&\RR_{23}(1,q)\RR_{13}(u,qv)\RR_{12}(u,v)|1,\ell,1\>=f(u^2,q^2v^2)\RR_{23}(1,q)|1,\ell,1\>+\\
&\quad+\frac{1}{4}\Big(\pf_{\ell1}(u,v)\pf_{1\ell}(u,qv)+\pf_{\ell1}(-u,v)\pf_{1\ell}(-u,qv)\Big)
\RR_{23}(1,q)|1,1,\ell\>+\\
&\quad+\frac14\Big(\pf_{\ell1}(u,v)\pf_{1\ell'}(-u,qv)+\pf_{\ell1}(-u,v)\pf_{1\ell'}(u,qv)\Big)
\RR_{23}(1,q)|1,1,\ell'\>
\end{split}
\end{equation}
Functions in front of $\RR_{23}(1,q)|1,1,\ell\>$ and $\RR_{23}(1,q)|1,1,\ell'\>$
in the second and third lines of this equality are equal to
\begin{equation*}
\frac{1+q}{2}\ \gri(u^2,v^2)\gle(u^2,q^2v^2)\quad\mbox{and}\quad
\frac{1-q}{2}\ \gri(u^2,v^2)\gle(u^2,q^2v^2)
\end{equation*}
respectively. Up to the multiplication of the second and third lines on the right hand side of
\r{em8} by the common  function $\gri(u^2,v^2)\gle(u^2,q^2v^2)$ these two lines are
equal to
\begin{equation}\label{em9}
\begin{split}
&\frac12\ \RR_{23}(1,q)\Big(|1,1,\ell\>+|1,1,\ell'\>\Big)+\frac{q}2\
\RR_{23}(1,q)\Big(|1,1,\ell\>-|1,1,\ell'\>\Big)=\\
&\quad=-\frac{q}2\ \RR_{23}(1,q)\Big(|1,\ell,1\>+|1,\ell',1\>\Big)
-\frac{q}2\ \RR_{23}(1,q)\Big(|1,\ell,1\>-|1,\ell',1\>\Big)=\\
&\qquad=-q\ \RR_{23}(1,q)|1,\ell,1\>
\end{split}
\end{equation}
where equalities
\begin{equation}\label{em99}
\RR(1,q)\Big(|1,\ell\>+|1,\ell'\>\Big)=-q \RR(1,q)\Big(|\ell,1\>+|\ell',1\>\Big)
\end{equation}
and
\begin{equation}\label{em98}
\RR(1,q)\Big(|1,\ell\>-|1,\ell'\>\Big)=- \RR(1,q)\Big(|\ell,1\>-|\ell',1\>\Big)
\end{equation}
have been used. The latter equalities are a consequence of the relations
\r{em4a}, \r{em4b} and \r{em4c}, \r{em4d} respectively.
Substituting \r{em9} to the equality \r{em8} we prove \r{em5a} for
$\ell=n,n+1$.
\qed

\begin{lemma}\label{5.3}
There are equalities for $1<i,j<2n$
\begin{subequations}\label{e8}
\begin{equation}\label{e8a}
\begin{split}
&\RR^{(n)}_{12}(1,q)\RR^{(n)}_{34}(1,q)\RR^{(n)}_{14}(u,qv)\RR^{(n)}_{13}(u,v)|i,1,j,1\>=\\
&\quad=\RR^{(n)}_{12}(1,q)\RR^{(n)}_{34}(1,q)\RR^{(n-1)}_{13}(u,v)|i,1,j,1\>
\end{split}
\end{equation}
and
\begin{equation}\label{e8b}
\begin{split}
&\<i,1,j,1|\RR^{(n)}_{13}(u,v)\RR^{(n)}_{14}(u,qv) \RR^{(n)}_{34}(1,q)\RR^{(n)}_{12}(1,q)=\\
&\quad=\<i,1,j,1|\RR^{(n-1)}_{13}(u,v)\RR^{(n)}_{34}(1,q)\RR^{(n)}_{12}(1,q)\,.
\end{split}
\end{equation}
\end{subequations}
\end{lemma}

Proof of this proposition is quite lengthy and can be found in Appendix~\ref{ApD}.\qed

\subsection{The embedding theorem}
\label{sec:embegging theorem}

\begin{thm}\label{main-st}
The commutation relations \r{RLL} in the algebra $U_q(\Dn{n})$ imply the
commutation relations in the algebra $U_q(\Dn{n-1})$ for the operators
$\MM^\pm(u)$ ($\mu,\rho=\pm$)
 \begin{equation}\label{RLL-red}
\RR^{(n-1)}(u,v) \cdot (\MM^\mu(u) \ot \Id)\cdot (\Id\ot \MM^\rho(v))=
(\Id\ot \MM^\rho(v))\cdot (\MM^\mu(u) \ot \Id)\cdot \RR^{(n-1)}(u,v)\,.
\end{equation}
\end{thm}

To prove theorem~\ref{main-st} we consider
$\RR\LL\LL$-commutation relation for fused $\LL$-operators $\Lf(u)$ and $\Lf(v)$ \eqref{Lf}
\begin{equation}\label{Em11}
\begin{split}
&\RR^{(n)}_{23}(qu,v)\RR^{(n)}_{13}(u,v)\RR^{(n)}_{24}(u,v)\RR^{(n)}_{14}(u,qv)\
\Lf^{(1,2)}(u)\ \Lf^{(3,4)}(v)=\\
&\quad =\Lf^{(3,4)}(v)\ \Lf^{(1,2)}(u)\ \RR^{(n)}_{14}(u,qv)\RR^{(n)}_{24}(u,v)
\RR^{(n)}_{13}(u,v)\RR^{(n)}_{23}(qu,v)\,.
\end{split}
\end{equation}
and   for $1<i_1,j_1,i_2,j_2<N$ take the matrix element of this commutation
relation
\begin{equation}\label{Em21}
\begin{split}
&\<i_1,1,j_1,1|
\RR^{(n)}_{23}(qu,v)\RR^{(n)}_{13}(u,v)\RR^{(n)}_{24}(u,v)\RR^{(n)}_{14}(u,qv)
\RR^{(n)}_{12}(1,q) \RR^{(n)}_{34}(1,q)\times\\
&\qquad \times \LL^{(1)}(u)\LL^{(2)}(qu) \LL^{(3)}(v)\LL^{(4)}(qv)
|i_2,1,j_2,1\>
=\\
&\quad =\<i_1,1,j_1,1|\LL^{(4)}(qv)\LL^{(3)}(v) \LL^{(2)}(qu)\LL^{(1)}(u)\times\\
&\qquad \times \RR^{(n)}_{34}(1,q) \RR^{(n)}_{12}(1,q)
 \RR^{(n)}_{14}(u,qv)\RR^{(n)}_{24}(u,v)\RR^{(n)}_{13}(u,v)\RR^{(n)}_{23}(qu,v)|i_2,1,j_2,1\>
\end{split}
\end{equation}
Let us transform last line in \r{Em21} using Proposition~\ref{5.3},
Lemma~\ref{5.2}
and Yang-Baxter
equation \r{YB}.
We have
\begin{equation}\label{Em22}
\begin{split}
&\RR^{(n)}_{34}(1,2) \RR^{(n)}_{12}(1,2)
 \RR^{(n)}_{14}(u,qv)\RR^{(n)}_{24}(u,v)\RR^{(n)}_{13}(u,v)\RR^{(n)}_{23}(qu,v)|i_2,1,j_2,1\>=\\
 &\quad= \RR^{(n)}_{12}(1,q)\RR^{(n)}_{13}(u,v) \RR^{(n)}_{14}(u,qv)
 \RR^{(n)}_{34}(1,q)\RR^{(n)}_{24}(u,v)\RR^{(n)}_{23}(qu,v)|i_2,1,j_2,1\>=\\
 &\quad=f(q^2u^2,v^2)\RR^{(n)}_{12}(1,q) \RR^{(n)}_{34}(1,q)
 \RR^{(n)}_{14}(u,qv)\RR^{(n)}_{13}(u,v)|i_2,1,j_2,1\>=\\
&\quad=f(q^2u^2,v^2)\RR^{(n)}_{12}(1,q) \RR^{(n)}_{34}(1,q)
\RR^{(n-1)}_{13}(u,v)|i_2,1,j_2,1\>\,.
\end{split}
\end{equation}
In the second step of this calculation we used equality \r{em5a} taken at $u\to qu$
and scaling invariance of $\RR$-matrix \r{scaling}.

Analogously first line in \r{Em21} can be transformed to
\begin{equation}\label{Em23}
\begin{split}
&\<i_1,1,j_1,1|
\RR^{(n)}_{23}(qu,v)\RR^{(n)}_{13}(u,v)\RR^{(n)}_{24}(u,v)\RR^{(n)}_{14}(u,qv)
\RR^{(n)}_{12}(1,q) \RR^{(n)}_{34}(1,q)=\\
&\quad = f(q^2u^2,v^2) \<i_1,1,j_1,1|
\RR^{(n-1)}_{13}(u,v)\RR^{(n)}_{12}(1,q) \RR^{(n)}_{34}(1,q)\,,
\end{split}
\end{equation}
where we used \r{em5b} at $v\to q^{-1}v$.

Equalities \r{Em22} and \r{Em23} allow to rewrite \r{Em21} in the form
\begin{equation*}\label{Em24}
\begin{split}
&\<i_1,1,j_1,1|\RR^{(n-1)}_{13}(u,v)\ \Lf^{(1,2)}(u)\Lf^{(3,4)}(v)|i_2,1,j_2,1\>=\\
&\qquad = \<i_1,1,j_1,1|\Lf^{(3,4)}(v)\Lf^{(1,2)}(u)\ \RR^{(n-1)}_{13}(u,v)|i_2,1,j_2,1\>
\end{split}
\end{equation*}
which proves the statement of  Theorem \r{RLL-red} due to Proposition~\ref{5.2}
and relation \r{e4}.

\qed

\section{Cartan-Weyl realization of $U_q(\Dn{n})$}
\label{cur-real}

In section~\ref{Algebra} the quantum loop algebra $U_q(\Dn{n})$
was formulated in terms of $\LL$-operators and $\RR\LL\LL$-relations.
But the matrix entries $\LL^\pm_{i,j}(u)$ \r{series} are not algebraically
independent.
Using methods of the papers
\cite{LiashykPakuliakRagoucySlavnov2019newsymmetry,LiashykPakuliak2022U_q(ABCD)}
and according to the relations \r{identalt} one can choose
the set of algebraically independent  generators in the quantum loop algebra $U_q(\Dn{n})$
as follows
\begin{equation}\label{gen-set}
\FF^\pm_{i+1,i}(u),\quad \EE^\pm_{i,i+1}(u),\quad i=1,\ldots,n-1,\quad
k^\pm_j(u),\quad j=1,\ldots,n\,.
\end{equation}
The central elements $z^\pm(u)$ \eqref{eq: zc z} in this case takes the form
\begin{equation}\label{cDn}
z^\pm(u)=k^\pm_n(u)k^\pm_n(-u)\prod_{\ell=1}^{n-1}
\frac{k^\pm_\ell(q^{n-\ell}u)}{k^\pm_\ell(q^{n-\ell-1}u)}\,.
\end{equation}

Since the modes of the Gaussian coordinates  $\FF^\pm_{i+1,i}(u)$,
$\EE^\pm_{i,i+1}(u)$
for $i=1,\ldots,n-1$ and $k^\pm_j(u)$, $j=1,\ldots,n$ are related to the
Cartan-Weyl generators of the quantum loop algebra $U_q(\Dn{n})$
\cite{KHOROSHKIN1993445} we call the realization of this algebra
in terms of the Gaussian coordinates \r{gen-set}
the {\it Cartan-Weyl realization}.

Below we will present
this realization  for small rank algebras
$U_q(\Dn{2})$ and $U_q(\Dn{3})$.
But before considering these small rank cases, one can look first
at the structure of monodromy matrix entries for the trivial algebra $U_q(\Dn{1})$.
$\RR$-matrix \r{tRalt} in this trivial case becomes diagonal
\begin{equation}\label{RD1}
\begin{split}
  \RR(u,v)  &= \fgo(u,v)\fgo(v,-u) \Big(\E_{11}\otimes \E_{11}+\E_{22}\otimes \E_{22} \Big)+ \\
 &\quad+ \fgo(u,-v)\fgo(v,u)(\E_{11}\otimes \E_{22} + \E_{22}\otimes \E_{11} ),
\end{split}
\end{equation}
where
\begin{equation}\label{fgo}
   \fgo(u,v)=\frac{q^{1/2}u-q^{-1/2}v}{u-v}\,.
\end{equation}
Then  $(2,1,1,1)$ element of the $\RR\LL\LL$ relations \r{RLL} has a form
\begin{equation*}
    \fgo(u,v)\fgo(v,-u) \LL_{1,1}(u) \LL_{2,1}(v)  =
    \fgo(u,-v)\fgo(v,u) \LL_{2,1}(v) \LL_{1,1}(u)\,.
\end{equation*}
At the point $u=qv$, this equality reduces to
\begin{equation*}
   \LL_{1,1}(qv) \LL_{2,1}(v) = 0\,,
\end{equation*}
which implies that $\EE_{1,2}(v) = 0$ due to the fact that diagonal monodromy
matrix element $\LL_{1,1}(u)$ is invertible. Analogously, one can check that
$\FF_{2,1}(u) = 0$. This proves that a trivial algebra $U_q(\Dn{1})$
is abelian and due to \r{Dpm} $\LL$-operator is a diagonal
matrix $\LL(u)={\rm diag}(k_1(u),k_1(-u))$. Since embedding theorem
is working also for this extremal case $n=2$, which means that
Gaussian coordinates $\EE_{2,3}(u)$ and $\FF_{3,2}(u)$ should vanish
in the algebra $U_q(\Dn{2})$. By induction and embedding theorem, one gets
 that $\EE_{n,n+1}(u)$ and $\FF_{n+1,n}(u)$ are also vanishing
in the algebra $U_q(\Dn{n})$.

\subsection{Algebra $U_q(\Dn{2})$}
\label{ApB}

For $n=2$ one has $\xi=q^{1-n}=q^{-1}$ and considering restrictions to
the matrix entries and Gaussian coordinates following from equalities
\r{D22-1} and \r{identalt} one concludes that algebraically independent
Gaussian coordinates in the algebra $U_q(\Dn{2})$ are $\EE^\pm_{1,2}(u)$,
$\FF^\pm_{2,1}(u)$, $k^\pm_1(u)$ and  $k^\pm_2(u)$
\begin{equation}\label{TF4}
 \mathbf{F}^\pm(u)=
\sk{\begin{array}{cccc}
1&\FF^\pm_{2,1}(u)&\FF^\pm_{2,1}(-u)&F^\pm_{4,1}(u)\\ [1mm]
0&1&0&-\FF^\pm_{2,1}(-u)\\[1mm] 0&0&1&-\FF^\pm_{2,1}(u)\\0&0&0&1
\end{array}}\,,
\end{equation}
\begin{equation}\label{TE4}
 \mathbf{E}^\pm(u)=
\sk{\begin{array}{cccc}
1&0&0&0\\ \EE^\pm_{1,2}(u)&1&0&0\\[1mm] \EE^\pm_{1,2}(-u)&0&1&0\\[1mm] \EE^\pm_{1,4}(u)
&-\EE^\pm_{1,2}(-u)&-\EE^\pm_{1,2}(u)&1
\end{array}}\,,
\end{equation}
and $\mathbf{K}^\pm(u)={\rm diag}(k^\pm_1(u),k^\pm_2(u),k^\pm_2(-u),k^\pm_4(u))$, where
$k^\pm_4(u)=z^\pm(u)k^\pm_1(qu)^{-1}$ and
\begin{equation}\label{TF5}
\begin{split}
\FF^\pm_{4,1}(u)&= \frac{q-1}{q+1}\ \FF^\pm_{2,1}(-u)^2 -\FF^\pm_{2,1}(-u)\FF^\pm_{2,1}(u)=\\
&= \frac{q-1}{q+1}\ \FF^\pm_{2,1}(u)^2 -\FF^\pm_{2,1}(u)\FF^\pm_{2,1}(-u)\,,
\end{split}
\end{equation}
\begin{equation}\label{TE5}
\begin{split}
\EE^\pm_{1,4}(u)&= \frac{1-q}{q+1}\ \EE^\pm_{1,2}(u)^2 -\EE^\pm_{1,2}(-u)\EE^\pm_{1,2}(u)=\\
&= \frac{1-q}{1+q}\ \EE^\pm_{1,2}(-u)^2 -\EE^\pm_{1,2}(u)\EE^\pm_{1,2}(-u)\,.
\end{split}
\end{equation}
Central elements $z^{\pm}(u)$ are equal to
\begin{equation*}
z^\pm(u)=k^\pm_2(u)k^\pm_2(-u)k^\pm_1(qu)k^\pm_1(u)^{-1}\,.
\end{equation*}

Befides rational function $\fgo(u,v)$ defined by the equality \r{fgo}
we introduce also functions
\begin{equation*}
\ggo(u,v)=\frac{(q^{1/2}-q^{-1/2})u}{u-v}\,,\quad
\tilde{\ggo}(u,v)=\frac{(q^{1/2}-q^{-1/2})v}{u-v}\,.
\end{equation*}
The commutation relations for the algebraically independent
generation series in the algebra $U_q(\Dn{2})$ can be obtained from the
$\RR\LL\LL$-commutation relations \r{RLL} and they are ($\mu,\nu=\pm$)
\begin{subequations}\label{kFE}
\begin{equation}\label{kFn1}
\begin{split}
k^\mu_{1}(v)\FF^\nu_{2,1}(u)k^\mu_{1}(v)^{-1}&=\fgo(u,v)\fgo(u,-v)\FF^\nu_{2,1}(u)-\\
&-
\frac12\Big(\gle(u,v)\FF^\mu_{2,1}(v)+\gle(u,-v)\FF^\mu_{2,1}(-v)\Big)\,,
\end{split}
\end{equation}
\begin{equation}\label{kFn2}
\begin{split}
k^\mu_{2}(v)\FF^\nu_{2,1}(u)k^\mu_{2}(v)^{-1}&=\fgo(v,u)\fgo(u,-v)\FF^\nu_{2,1}(u)+\\
&+\frac12\Big(\gle(u,v)\FF^\mu_{2,1}(v)-\gle(u,-v)\FF^\mu_{2,1}(-v)\Big)\,,
\end{split}
\end{equation}
\begin{equation}\label{kEn1}
\begin{split}
k^\mu_{1}(v)^{-1}\EE^\nu_{1,2}(u)k^\mu_{1}(v)&=\fgo(u,v)\fgo(u,-v)\EE^\nu_{1,2}(u)-\\
&-\frac12\Big(\gri(u,v)\EE^\mu_{1,2}(v)+\gri(u,-v)\EE^\mu_{1,2}(-v)\Big)\,,
\end{split}
\end{equation}
\begin{equation}\label{kEn2}
\begin{split}
k^\mu_{2}(v)^{-1}\EE^\nu_{1,2}(u)k^\mu_{2}(v)&=\fgo(v,u)\fgo(u,-v)\EE^\nu_{1,2}(u)+\\
&+\frac12\Big(\gri(u,v)\EE^\mu_{1,2}(v)-\gri(u,-v)\EE^\mu_{1,2}(-v)\Big)\,,
\end{split}
\end{equation}
\begin{equation}\label{En-Fn}
[\EE^\mu_{1,2}(u),\FF^\nu_{2,1}(v)]=\frac{\gle(v,u)}{2}\sk{\frac{k^\nu_{2}(v)}{k^\nu_{1}(v)}-
\frac{k^\mu_{2}(u)}{k^\mu_{1}(u)}}\,,
\end{equation}
\begin{equation}\label{FnFn}
\fgo(v,u)\FF^\mu_{2,1}(u)\FF^\nu_{2,1}(v)=\fgo(u,v) \FF^\nu_{2,1}(v)\FF^\mu_{2,1}(u)
+\ggo(v,u) \FF^\mu_{2,1}(u)^2+\tilde{\ggo}(v,u) \FF^\nu_{2,1}(v)^2\,,
\end{equation}
\begin{equation}\label{EnEn}
\fgo(u,v)\EE^\mu_{1,2}(u)\EE^\nu_{1,2}(v)=\fgo(v,u) \EE^\nu_{1,2}(v)\EE^\mu_{1,2}(u)
+\ggo(u,v) \EE^\mu_{1,2}(u)^2+\tilde{\ggo}(u,v) \EE^\nu_{1,2}(v)^2\,.
\end{equation}
\end{subequations}

\begin{prop}\label{isoD2A1}
Quantum loop algebra $U_q(\Dn{2})$ given by the commutation relations \r{kFE}
is isomorphic to the quantum loop algebra $U_{q^{1/2}}(A^{(1)}_1)$ with following
identification of the Gaussian coordinates
\begin{equation*}
\begin{split}
k^\pm_{1}(u|D^{(2)}_2)&=k^\pm_{1}(u|A^{(1)}_1)k^\pm_{1}(-u|A^{(1)}_1)\,,\\
k^\pm_{2}(u|D^{(2)}_2)&=k^\pm_{2}(u|A^{(1)}_1)k^\pm_{1}(-u|A^{(1)}_1)\,,\\
\EE^\pm_{1,2}(u|D^{(2)}_2)&=\Big(\frac{q^{1/2}+q^{-1/2}}{2}\Big)^{1/2}
\EE^\pm_{1,2}(u|A^{(1)}_1)\,,\\
\FF^\pm_{2,1}(u|D^{(2)}_2)&=\Big(\frac{q^{1/2}+q^{-1/2}}{2}\Big)^{1/2}
\FF^\pm_{2,1}(u|A^{(1)}_1)\,.
\end{split}
\end{equation*}
\end{prop}
Proof is by the direct verification using commutation relations between Gaussian coordinates
in the algebra $U_{q^{1/2}}(A^{(1)}_1)$.\qed

\subsection{Algebra $U_q(\Dn{3})$}
\label{ApC}

For $n=3$
\begin{equation*}
\Der={\rm diag}(q^{3/2},q^{1/2},1,1,q^{-1/2},q^{-3/2}),\quad \xi=q^{-2}
\end{equation*}
and equality \r{identalt} yields the relation between matrix entries of $\LL$-operators
$\LL(u)$ and $\hLL(u)$ of the form
\begin{equation}\label{D3-8}
q^{\bar\imath-\bar\jmath}\ \hLL^\pm_{i,j}(u)=z(u)^{-1}\LL^\pm_{i,j}(uq^2).
\end{equation}
Using this equality and \r{D22-1} one can find that
algebraically independent set of the generating series in quantum
loop algebra $U_q(\Dn{3})$ is given by \r{gen-set} for $n=3$.

Using $\RR\LL\LL$ commutation relations \r{RLL} one can find that
\begin{equation}\label{D3d1}
\mathbf{K}^\pm(u)
={\rm diag}\Big(k^\pm_{1}(u),k^\pm_2(u),k^\pm_3(u),k^\pm_3(-u),k^\pm_5(u),k^\pm_6(u)\Big)\,,
\end{equation}
where
\begin{equation}\label{D3-9}
k^\pm_6(u)=z^\pm(u)k^\pm_1(q^2u)^{-1}
\end{equation}
and
\begin{equation}\label{D3-12}
k^\pm_5(u)=z^\pm(u)\ \frac{1}{k^\pm_2(qu)}\ \frac{k^\pm_1(qu)}{k^\pm_1(q^2u)}\,.
\end{equation}
The central elements $z^\pm(u)$ are given by equality \r{cDn} for $n=3$.

Upper- and lower-triangular matrices $\mathbf{F}^\pm(u)$ and $\mathbf{E}^\pm(u)$ are
\begin{equation}\label{D3d3}
\mathbf{F}^\pm(u)=
\sk{\begin{array}{cccccc}
1&\FF^\pm_{2,1}(u)&\FF^\pm_{3,1}(u)&\FF^\pm_{3,1}(-u)
&\FF^\pm_{5,1}(u)&\FF^\pm_{6,1}(u)\\[1mm]
0&1&\FF^\pm_{3,2}(u)&\FF^\pm_{3,2}(-u)&\FF^\pm_{5,2}(u)&\FF^\pm_{6,2}(u)\\[1mm]
0&0&1&0&-\FF^\pm_{3,2}(-u)&\FF^\pm_{6,4}(-u)\\[1mm]
0&0&0&1&-\FF^\pm_{3,2}(u)&\FF^\pm_{6,4}(u)\\[1mm]
0&0&0&0&1&\FF^\pm_{6,5}(u)\\[1mm]
0&0&0&0&0&1
\end{array}}
\end{equation}
\begin{equation}\label{D3d4}
\mathbf{E}^\pm(u)=
\sk{\begin{array}{cccccc}
1&0&0&0&0&0\\[1mm]
\EE^\pm_{1,2}(u) &1&0&0&0&0\\[1mm]
\EE^\pm_{1,3}(u) &\EE^\pm_{2,3}(u)&1&0&0&0\\[1mm]
\EE^\pm_{1,3}(-u)&\EE^\pm_{2,3}(-u)&0&1&0&0\\[1mm]
\EE^\pm_{1,5}(u)&\EE^\pm_{2,5}(u)&-\EE^\pm_{2,3}(-u)&-\EE^\pm_{2,3}(u)&1&0\\[1mm]
\EE^\pm_{1,6}(u)&\EE^\pm_{2,6}(u)&\EE^\pm_{4,6}(-u)&\EE^\pm_{4,6}(u)&\EE^\pm_{5,6}(u)&1
\end{array}}\,,
\end{equation}
where
\begin{equation}\label{D3d11}
\begin{split}
\FF^\pm_{5,2}&=\frac{q-1}{q+1}\FF^\pm_{3,2}(u)^2-\FF^\pm_{3,2}(u)\FF^\pm_{3,2}(-u)\,,\\
\EE^\pm_{2,5}&=\frac{1-q}{q+1}\EE^\pm_{2,3}(-u)^2-\EE^\pm_{2,3}(u)\EE^\pm_{2,3}(-u)\,.
\end{split}
\end{equation}
One can note that matrix elements $\FF_{5,2}(u)$, $\EE_{2,5}(u)$
coincide with the elements $\FF_{4,1}(u)$, $\EE_{1,4}(u)$
calculated in case of $n=2$ by the formulas \r{TF5} and \r{TE5}.
This demonstrates an embedding of $U_q(\Dn{2})$ into $U_q(\Dn{3})$.
One can also got expressions for other Gaussian coordinates in \r{D3d3} and
\r{D3d4} through independent set \r{gen-set} for $n=3$, but this is not important
for further discussion.

Below we list only nontrivial commutation relations between Gaussian coordinates \r{gen-set}
for $n=3$
in the quantum loop algebra $U_q(\Dn{3})$ ($\mu,\nu=\pm$)
\begin{subequations}\label{D3cr}
\begin{equation}\label{k12F}
\begin{split}
k^\mu_{1}(v)\FF^\nu_{2,1}(u)k^\mu_{1}(v)^{-1}&=
f(u^2,v^2)\FF^\nu_{2,1}(u)-\gle(u^2,v^2)\FF^\nu_{2,1}(v)\,,\\
k^\mu_{2}(v)\FF^\nu_{2,1}(u)k^\mu_{2}(v)^{-1}&=
f(v^2,u^2)\FF^\nu_{2,1}(u)-\gri(v^2,u^2)\FF^\nu_{2,1}(v)\,,
\end{split}
\end{equation}

\begin{equation}\label{k23F}
\begin{split}
k^\mu_{2}(v)\FF^\nu_{3,2}(u)k^\mu_{2}(v)^{-1}&=\fgo(u,v)\fgo(u,-v)\FF^\nu_{3,2}(u)-\\
&-\frac12\Big(\gle(u,v)\FF^\mu_{3,2}(v)+\gle(u,-v)\FF^\mu_{3,2}(-v)\Big)\,,\\
k^\mu_{3}(v)\FF^\nu_{3,2}(u)k^\mu_{3}(v)^{-1}&=\fgo(v,u)\fgo(u,-v)\FF^\nu_{3,2}(u)+\\
&+\frac12\Big(\gle(u,v)\FF^\mu_{3,2}(v)-\gle(u,-v)\FF^\mu_{3,2}(-v)\Big)\,,
\end{split}
\end{equation}

\begin{equation}\label{k12E}
\begin{split}
k^\mu_{1}(v)^{-1}\EE^\nu_{1,2}(u)k^\mu_{1}(v)&=
f(u^2,v^2)\EE^\nu_{1,2}(u)-\gri(u^2,v^2)\EE^\nu_{1,2}(v)\,,\\
k^\mu_{2}(v)^{-1}\EE^\nu_{1,2}(u)k^\mu_{2}(v)&=
f(v^2,u^2)\EE^\nu_{1,2}(u)-\gle(v^2,u^2)\EE^\nu_{1,2}(v)\,,
\end{split}
\end{equation}

\begin{equation}\label{k23E}
\begin{split}
k^\mu_{2}(v)^{-1}\EE^\nu_{2,3}(u)k^\mu_{2}(v)&=\fgo(u,v)\fgo(u,-v)\EE^\nu_{2,3}(u)-\\
&-\frac12\Big(\gri(u,v)\EE^\mu_{2,3}(v)+\gri(u,-v)\EE^\mu_{2,3}(-v)\Big)\,,\\
k^\mu_{3}(v)^{-1}\EE^\nu_{2,3}(u)k^\mu_{3}(v)&=\fgo(v,u)\fgo(u,-v)\EE^\nu_{2,3}(u)+\\
&+\frac12\Big(\gri(u,v)\EE^\mu_{2,3}(v)-\gri(u,-v)\EE^\mu_{2,3}(-v)\Big)\,,
\end{split}
\end{equation}

\begin{equation}\label{3EF}
\begin{split}
[\EE^\mu_{1,2}(u),\FF^\nu_{2,1}(v)]&=\gle(v^2,u^2)\sk{\frac{k^\nu_{2}(v)}{k^\nu_{1}(v)}-
\frac{k^\mu_{2}(u)}{k^\mu_{1}(u)}}\,,\\
[\EE^\mu_{2,3}(u),\FF^\nu_{3,2}(v)]&=\frac{\gle(v,u)}{2}\sk{\frac{k^\nu_{3}(v)}{k^\nu_{2}(v)}-
\frac{k^\mu_{3}(u)}{k^\mu_{2}(u)}}\,,
\end{split}
\end{equation}

\begin{equation}\label{3FF}
\begin{split}
f(v^2,u^2)\FF^\mu_{2,1}(u)\FF^\nu_{2,1}(v)&=f(u^2,v^2) \FF^\nu_{2,1}(v)\FF^\mu_{2,1}(u)+\\
&+\gle(v^2,u^2) \FF^\mu_{2,1}(u)^2+\gri(v^2,u^2) \FF^\nu_{2,1}(v)^2\,,\\
\fgo(v,u)\FF^\mu_{3,2}(u)\FF^\nu_{3,2}(v)&=\fgo(u,v) \FF^\nu_{3,2}(v)\FF^\mu_{3,2}(u)+\\
&+\ggo(v,u) \FF^\mu_{3,2}(u)^2+\tilde{\ggo}(v,u) \FF^\nu_{3,2}(v)^2\,,\\
\FF^\mu_{2,1}(u)\FF^\nu_{3,2}(v)&=f(v^2,u^2)\FF^\nu_{3,2}(v)\FF^\mu_{2,1}(u)-
\gri(v^2,u^2)\FF^\nu_{3,2}(v)\FF^\nu_{2,1}(v)+\\
+\gri(v^2,u^2)&\FF^\nu_{3,1}(v)+\frac12\Big(\gri(u,v)\FF^\mu_{3,1}(u)+
\gri(u,-v)\FF^\mu_{3,1}(-u)\Big)\,,
\end{split}
\end{equation}

\begin{equation}\label{3EE}
\begin{split}
f(u^2,v^2)\EE^\mu_{1,2}(u)\EE^\nu_{1,2}(v)&=f(v^2,u^2) \EE^\nu_{1,2}(v)\EE^\mu_{1,2}(u)+\\
&+\gle(u^2,v^2) \EE^\mu_{1,2}(u)^2+\gri(u^2,v^2) \EE^\nu_{1,2}(v)^2\,,\\
\fgo(u,v)\EE^\mu_{2,3}(u)\EE^\nu_{2,3}(v)&=\fgo(v,u) \EE^\nu_{2,3}(v)\EE^\mu_{2,3}(u)+\\
&+\ggo(u,v) \EE^\mu_{2,3}(u)^2+\tilde{\ggo}(u,v) \EE^\nu_{2,3}(v)^2\,,\\
\EE^\mu_{2,3}(u)\EE^\nu_{1,2}(v)&=f(u^2,v^2)\EE^\nu_{1,2}(v)\EE^\mu_{2,3}(u)-
\gle(u^2,v^2)\EE^\mu_{1,2}(u)\EE^\mu_{2,3}(u)+\\
+\gle(u^2,v^2)&\EE^\mu_{1,3}(u)+\frac12\Big(\gle(v,u)\EE^\nu_{1,3}(v)+
\gri(v,-u)\EE^\nu_{1,3}(-v)\Big)\,.
\end{split}
\end{equation}
\end{subequations}

\subsection{Isomorphism of $U_q(\Dn{3})$ and $U_{q^{1/2}}(A^{(2)}_3)$}

As well as quantum loop algebra $U_q(\Dn{2})$ the algebra $U_q(\Dn{3})$
is isomorphic to the algebra of another series.

\begin{prop}\label{isoD3A3}
Quantum loop algebra $U_q(\Dn{3})$ given by the commutation relations \r{D3cr}
is isomorphic to the quantum loop algebra $U_{q^{1/2}}(A^{(2)}_3)$ with following
identification of the Gaussian coordinates
\begin{equation}
\begin{split}
k^\pm_{1}(u|D^{(2)}_3)&=k^\pm_{1}(u|A^{(2)}_3)k^\pm_{1}(-u|A^{(2)}_3)
k^\pm_{2}(q^{-1}u|A^{(2)}_3)\,,\\\
k^\pm_{2}(u|D^{(2)}_3)&=k^\pm_{1}(u|A^{(2)}_3)k^\pm_{1}(-u|A^{(2)}_3)
k^\pm_{3}(q^{-1}u|A^{(2)}_3)\,,\\
k^\pm_{3}(u|D^{(2)}_3)&=k^\pm_{2}(u|A^{(2)}_3)k^\pm_{1}(-u|A^{(2)}_3)
k^\pm_{3}(q^{-1}u|A^{(2)}_3)\,,\\
\EE^\pm_{1,2}(u|D^{(2)}_3)&=u^{-1}\EE^\pm_{2,3}(q^{-1}u|A^{(2)}_3)\,,\\
\FF^\pm_{2,1}(u|D^{(2)}_3)&=u\ \FF^\pm_{3,2}(q^{-1}u|A^{(2)}_3)\,,\\
\EE^\pm_{2,3}(u|D^{(2)}_3)&=\Big(\frac{q^{1/2}+q^{-1/2}}{2}\Big)^{1/2}
\EE^\pm_{1,2}(u|A^{(2)}_3)\,,\\
\FF^\pm_{2,3}(u|D^{(2)}_3)&=\Big(\frac{q^{1/2}+q^{-1/2}}{2}\Big)^{1/2}
\FF^\pm_{2,1}(u|A^{(2)}_3)\,.
\end{split}
\end{equation}
\end{prop}
The proof is by the direct verification using commutation relations between Gaussian coordinates
in the algebra $U_{q^{1/2}}(A^{(2)}_3)$ which can be found in
\cite{LiashykPakuliak2022U_q(ABCD)}.\qed

\subsection{Algebra $U_q(\Dn{n})$}

Algebraically independent set of Gaussian coordinates which generate the algebra
$U_q(\Dn{n})$ is given by \r{gen-set} with
 the central elements \r{cDn}.
Corresponding to the simple roots
dependent Gaussian coordinates  can be obtained from
equalities \r{D22-1} and \r{identalt} and are
\begin{equation}\label{depEF}
\EE^\pm_{(i+1)',i'}(u)=-\EE^\pm_{i,i+1}(u\xi^{-1}q^{-i}),\quad
\FF^\pm_{i',(i+1)'}(u)=-\FF^\pm_{i+1,i}(u\xi^{-1}q^{-i})
\end{equation}
for $i=1,\ldots,n-1$ and
\begin{equation}\label{depk}
k^\pm_{i'}(u)=\frac{k^\pm_n(u)k^\pm_n(-u)}{k^\pm_i(u\xi^{-1}q^{1-i})}
\prod_{\ell=i}^{n-1}\frac{k^\pm_\ell(u\xi^{-1}q^{1-\ell})}{k^\pm_\ell(u\xi^{-1}q^{-\ell})}
\end{equation}
for $i=1,\ldots,n$. Note that
\r{depk} for $i=n$ yields an equality $k^\pm_{n+1}(u)=k^\pm_n(-u)$.

Following the Ding-Frenkel approach \cite{DingFrenkel1993isomorphism}
and  the embedding theorem~\ref{main-st}
one can write the commutation relations in the algebra $U_q(\Dn{n})$
using the {\it currents}
\begin{equation}\label{DF}
E_i(u)=\EE^+_{i,i+1}(u)-\EE^-_{i,i+1}(u)\,,\quad
F_i(u)=\FF^+_{i+1,i}(u)-\FF^-_{i+1,i}(u)
\end{equation}
for $i=1,\ldots,n-1$ and the Cartan currents $k^\pm_j(u)$ for $j=1,\ldots,n$.

 The currents
$E_i(u)$, $F_i(u)$, $i=1,\ldots,n-2$ and $k^\pm_j(u)$, $j=1,\ldots,n-1$ are
series with respect to the square of the spectral parameter $u^2$. They have
$U_q(\mathfrak{gl}_{n-1})$-type commutation relations.
We formulate  all nontrivial commutation relations between these currents without proofs
\begin{equation*}\label{kiFA}
\begin{split}
k^{\pm}_i(u) F_i(v) k^{\pm}_i(u)^{-1}&= \frac{q^{-1}u^2 -qv^2}{u^2-v^2}\   F_i(v),\\
k^{\pm}_{i+1}(u)F_i(v)k^{\pm}_{i+1}(u)^{-1}&= \frac{qu^2 -q^{-1}v^2}{u^2-v^2}\   F_i(v),
\end{split}
\end{equation*}
\begin{equation}\label{kEFA}
\begin{split}
k^{\pm}_i(u)^{-1}E_i(v)k^{\pm}_i(u)&=\frac{q^{-1}u^2 -qv^2}{u^2-v^2}\   E_i(v),\\
k^{\pm}_{i+1}(u)^{-1}E_i(v)k^{\pm}_{i+1}(u)&=\frac{qu^2 -q^{-1}v^2}{u^2-v^2}\  E_i(v),
\end{split}
\end{equation}
\begin{equation*}\label{FiFiA}
f(v^2,u^2)\ F_i(u)F_i(v)=  f(u^2,v^2)\  F_i(v)F_i(u),
\end{equation*}
\begin{equation*}\label{EiEiA}
f(u^2,v^2)\ E_i(u) E_i(v)=  f(v^2,u^2)\  E_i(v) E_i(u),
\end{equation*}
\begin{equation*}\label{FiFiiA}
(u^2-v^2)\ F_i(u)F_{i+1}(v)= (q^{-1}u^2-qv^2)\ F_{i+1}(v)F_i(u),
\end{equation*}
\begin{equation*}\label{EiEiiA}
(q^{-1}u^2-qv^2)\ E_i(u)E_{i+1}(v)= (u^2-v^2)\  E_{i+1}(v)E_i(u),
\end{equation*}
\begin{equation*}\label{EFA}
[E_i(u),F_j(v)]=\delta_{i,j}\ (q-q^{-1})
\delta(u^2,v^2)\Big(k^-_{i+1}(v)\,k^-_{i}(v)^{-1}-k^+_{i+1}(u)\,k^{+}_{i}(u)^{-1}\Big).
\end{equation*}
In last equality $\delta(x,y)$ is a multiplicative $\delta$-function given by a series
\begin{equation}\label{delta}
\delta(x,y)=\sum_{i\in\ZZ} x^i y^{-i}\,.
\end{equation}
 There are also  Serre relations for the currents $E_i(u)$ and $F_i(u)$
 \cite{Drinfeld1988NewRealisation,DingFrenkel1993isomorphism}
\begin{equation}\label{Serre}
\begin{split}
\mathop{\rm Sym}_{v_1,v_2}\Big[F_i(v_1),[F_i(v_2),F_{i\pm1}(u)]_{q^{-1}}\Big]_q&=0\,,\\
\mathop{\rm Sym}_{v_1,v_2}\Big[E_i(v_1),[E_i(v_2),E_{i\pm1}(u)]_{q}\Big]_{q^{-1}}&=0\,,
\end{split}
\end{equation}
where  $[A,B]_q$ means $q$-commutator
\begin{equation*}
[A,B]_q=A\,B-q\ B\,A
\end{equation*}
and $\mathop{\rm Sym}_{v_1,v_2}G(v_1,v_2)\equiv G(v_1,v_2)+G(v_2,v_1)$.

The currents $E_{n-1}(u)$, $F_{n-1}(u)$ and $k_n(u)$ are currents with respect to the
spectral parameter $u$ and nontrivial commutation relations between them
and other currents are
\begin{equation}\label{kEFD}
\begin{split}
k^\pm_{n-1}(u)F_{n-1}(v)k^\pm_{n-1}(u)^{-1}&=\frac{(u-q^{-1}v)(v+qu)}{u^2-v^2}F_{n-1}(v)\,,\\
k^\pm_{n}(u)F_{n-1}(v)k^\pm_{n}(u)^{-1}&=\frac{(u-qv)(u+q^{-1}v)}{u^2-v^2}F_{n-1}(v)\,,\\
k^\pm_{n-1}(u)^{-1}E_{n-1}(v)k^\pm_{n-1}(u)&=\frac{(u-q^{-1}v)(v+qu)}{u^2-v^2}E_{n-1}(v)\,,\\
k^\pm_{n}(u)^{-1}E_{n-1}(v)k^\pm_{n}(u)&=\frac{(u-qv)(u+q^{-1}v)}{u^2-v^2}E_{n-1}(v)\,,\\
\fgo(v,u)\ F_{n-1}(u)F_{n-1}(v)&=\fgo(u,v)\ F_{n-1}(v)F_{n-1}(u)\,,\\
\fgo(u,v)\ E_{n-1}(u)E_{n-1}(v)&=\fgo(v,u)\ E_{n-1}(v)E_{n-1}(u)\,,\\
(u^2-v^2)\ F_{n-2}(u)F_{n-1}(v)&= (q^{-1}u^2-qv^2)\ F_{n-1}(v)F_{n-2}(u)\,,\\
(q^{-1}u^2-qv^2)\ E_{n-2}(u)E_{n-1}(v)&= (u^2-v^2)\  E_{n-1}(v)E_{n-2}(u)\,,\\
[E_{n-1}(u),F_{n-1}(v)]=\frac12(q-q^{-1})&\delta(u,v)
\Big(k^-_{n}(v)\,k^-_{n-1}(v)^{-1}-k^+_{n}(u)\,k^{+}_{n-1}(u)^{-1}\Big)
\end{split}
\end{equation}
and the Serre relations
\begin{equation}\label{SerreDn}
\begin{split}
\mathop{\rm Sym}_{v_1,v_2} \Big[F_{n-2}(v_1),
[F_{n-2}(v_2),F_{n-1}(u)]_{q^{-1}}\Big]_{q} &=0\,,\\
\mathop{\rm Sym}_{v_1,v_2}
\Big[E_{n-2}(v_1),[E_{n-2}(v_2),E_{n-1}(u)]_{q}\Big]_{q^{-1}} &=0\,,\\
\mathop{\rm Sym}_{v_1,v_2,v_3}\Big[F_{n-1}(v_1),\Big[F_{n-1}(v_2),[F_{n-1}(v_3),F_{n-2}(u)]_{q^{-1}}\Big]_{q}\Big]&=0\,,\\
\mathop{\rm Sym}_{v_1,v_2,v_3}\Big[E_{n-1}(v_1),\Big[E_{n-1}(v_2),[E_{n-1}(v_3),E_{n-2}(u)]_{q}\Big]_{q^{-1}}\Big]&=0\,.
\end{split}
\end{equation}

Note that according the commutation relations \r{kEFD} the currents $F_{n-1}(u)$ and
$F_{n-1}(-u)$ commute (analogously, $[E_{n-1}(u),E_{n-1}(-u)]=0$).
But this does not mean that the Gaussian coordinates $\FF^\mu_{n,n-1}(u)$ and
$\FF^\nu_{n,n-1}(-u)$ are commutative for $\mu,\nu=\pm$.
In particular, this non-commutativity provides the equalities in \r{TF5} and \r{TE5}.


\subsection{Composed currents and Gaussian coordinates}

Denote by  $\overline U_f$
an extension of the algebra $U_f=U^-_f\,\cup\, U^+_f\,\cup\, U^+_k$ formed
by linear combinations of series, given as infinite sums of monomials
$a_{i_1}[n_1]\cdots a_{i_k}[n_k]$ with $n_1\leq\cdots\leq n_k$, and $n_1+...+n_k$
fixed, where  $a_{i_l}[n_l]$ is either $F_{i_l}[n_l]$ or $k^+_{i_l}[n_l]$.
Analogously, denoted by
$\overline U_e$ an extension of the algebra $U_e=U^-_e\,\cup\, U^+_e\,\cup\, U^-_k$ formed
by linear combinations of series, given as infinite sums of monomials
$a_{i_1}[n_1]\cdots a_{i_k}[n_k]$ with $n_1\geq\cdots\geq n_k$, and $n_1+...+n_k$
fixed, where  $a_{i_l}[n_l]$ is either $E_{i_l}[n_l]$ or $k^-_{i_l}[n_l]$.

For $n<i<2n$ we introduce dependent total currents using formulas \r{DF} and
\r{depEF}
\begin{equation}\label{depcu}
\begin{split}
E_i(u)&=-E_{(i+1)'}(u\xi^{-1}q^{-(i+1)'})=-E_{(i+1)'}(uq^{i-n-1})\,,\\
F_i(u)&=-F_{(i+1)'}(u\xi^{-1}q^{-(i+1)'})=-F_{(i+1)'}(uq^{i-n-1})\,.
\end{split}
\end{equation}

Let us introduce so called {\it composed currents} $\Fn_{j,i}(u)$ and
$\En_{i,j}(u)$ for $1\leq i<j\leq 2n$
\begin{equation}\label{comp-cuF}
\Fn_{j,i}(u)=\begin{cases}
F_{j-1}(u)F_{j-2}(u)\cdots F_{i+1}(u)F_{i}(u),&
i<j< n+1,\quad n< i<j\,,\\
F_{n-1}(-u)F_{n-2}(u)\cdots F_i(u),&i< n,\quad j=n+1\,,\\
0,& i=n,\quad j=n+1\,,\\
F_{j-1}(u)\cdots F_{n+2}(u)F_{n+1}(-u),& i=n,\quad n+1< j\,,\\
F_{j-1}(u)\cdots F_{n+2}(u)F_{n+1}(-u)\times\\
\qquad\times F_{n-1}(u)F_{n-2}(u)\cdots F_i(u),&
i< n,\quad n+1<j\,,
\end{cases}
\end{equation}
\begin{equation}\label{comp-cuE}
\En_{i,j}(u)=\begin{cases}
E_{i}(u)E_{i+1}(u)\cdots E_{j-2}(u)E_{j-1}(u),&
i<j< n+1,\quad n< i<j\,,\\
E_{i}(u)\cdots E_{n-2}(u)E_{n-1}(-u),&i<n,\quad j=n+1\,,\\
0,& i=n,\quad j=n+1\,,\\
E_{n+1}(-u)E_{n+2}(u)\cdots E_{j-1}(u),& i=n,\quad n+1< j\,,\\
E_{i}(u)\cdots E_{n-2}(u)E_{n-1}(u)\times\\
\qquad\times E_{n+1}(-u)E_{n+2}(u)\cdots E_{j-1}(u),&
i< n,\quad n+1<j\,.
\end{cases}
\end{equation}

The composed currents are well-defined in the category of the highest weight
representations and are elements from completions $\overline{U}_e$ and
$\overline{U}_f$.
One can prove that
the actions of the projections
$\Pfpm$ and $\Pepm$ investigated in \cite{EnriquezKhoroshkinPakuliak2007Projections}
and shortly discussed in the section~\ref{nor-ord}
can be extended to completions $\overline{U}_f$ and
$\overline{U}_e$. These projections are well-defined when applied  to the composed currents
\r{comp-cuF} and \r{comp-cuE} respectively. Definition of the projections
depends on the ordering of the subalgebras formed by the Gaussian coordinates
of the same type \cite{EnriquezKhoroshkinPakuliak2007Projections}
and here we will use projections defined for the ordering \r{order1}.

In order to apply these projections to the products of the simple roots
currents, one has to do the following.
First, one should replace all currents by Ding-Frenkel linear combinations
\r{DF} and then order each monomial according to the ordering \r{order1}
using the commutation relations between them. After that, an application
of the projection $\Pfp$ to the normal ordered product of currents
is removing all the terms which have at least one 'negative' Gaussian coordinate
$\FF^-_{j,i}(u)$ on the left. Analogously, an application of the projection $\Pfm$
to the same normal ordering product of the currents is removing
all the terms which have at least one 'positive' Gaussian coordinate
$\FF^+_{j,i}(u)$ on the right.

Composed currents are related to the Gaussian coordinates as formulated in the following
proposition.
\begin{prop}\label{ccvsGC}
Gaussian coordinates can be related to the composed current through projections
$\Pfpm$ and $\Pepm$
\begin{equation}\label{gccc}
\begin{split}
\Pfp\big(\Fn_{j,i}(u)\big)&=\FF^+_{j,i}(u),\qquad
\Pfm\big(\Fn_{j,i}(u)\big)={\tFF}^-_{j,i}(u)\,,\\
\Pep\big(\En_{i,j}(u)\big)&=\EE^+_{i,j}(u),\qquad
\Pem\big(\En_{i,j}(u)\big)={\tEE}^-_{i,j}(u)\,.
\end{split}
\end{equation}
\end{prop}

Proof can be done along the same lines as in \cite{LiashykPakuliak2022U_q(ABCD)}.
Instead of repeating this proof, we give an example
and calculate $\Pfp\big(\Fn_{n+2,n-1}(u)\big)=\Pfp\big(F_{n+1}(-u)F_{n-1}(u)\big)$ which
correspond to the last line in \r{comp-cuF}
\begin{equation*}
\begin{split}
\Pfp\Big(F_{n+1}(-u)F_{n-1}(u)\Big)&=\Pfp\Big(\FF^+_{n,n-1}(-u)\FF^-_{n,n-1}(u)\Big)-
\FF^+_{n,n-1}(-u)\FF^+_{n,n-1}(u)=\\
&=\frac{\ggo(u,-u)}{\fgo(u,-u)}\FF^+_{n,n-1}(-u)^2-\FF^+_{n,n-1}(-u)\FF^+_{n,n-1}(u)=\\
&=\frac{q-1}{q+1}\FF^+_{n,n-1}(-u)^2-\FF^+_{n,n-1}(-u)\FF^+_{n,n-1}(u)=
\FF^+_{n+2,n-1}(u)
\end{split}
\end{equation*}
which is in accordance with the formula \r{TF5}.

In order to calculate the projection $\Pfp\Big(\FF^+_{n,n-1}(-u)\FF^-_{n,n-1}(u)\Big)$
we used analog of the relations \r{FnFn} and  \r{3FF} for the algebra $U_q(\Dn{n})$
\begin{equation*}
\begin{split}
\fgo(v,u)\FF^+_{n,n-1}(u)\FF^-_{n,n-1}(v)&=\fgo(u,v) \FF^-_{n,n-1}(v)\FF^+_{n,n-1}(u)+\\
&+\ggo(v,u) \FF^+_{n,n-1}(u)^2+\tilde{\ggo}(v,u) \FF^-_{n,n-1}(v)^2
\end{split}
\end{equation*}
specialized to the values of the spectral parameters $v=-u$.

Note, that products of the currents $F_{n+1}(-u)F_{n-1}(u)$ is formal series
with respect to spectral parameter $u^2$ since it is equal to
$F_{n+1}(u)F_{n-1}(-u)$ due to the commutativity of the currents
$F_{n-1}(u)$ and $F_{n-1}(-u)$ (recall that by definition \r{depcu}
$F_{n+1}(u)=-F_{n-1}(u)$).

Moreover, one can verify from the equalities \r{gccc} that Gaussian coordinates
$\FF^\pm_{j,i}(u)$ and $\EE^\pm_{i,j}(u)$
for $i<j<n$, $n<i<j$ and $i<n<j-1$ are series with respect
to $u^{\mp2}$. Together with the fact that series $k^\pm_j(u)$ for
$j<n$ and $n+1<j$
are also series with respect
to $u^{\mp2}$ we may conclude that matrix entries $\LL^\pm_{i,j}(u)$
have the same analytical properties for either $i\not=n,n+1$ or
$j\not=n,n+1$ and depend on the square of the spectral parameter.
This fact is in correspondence with the relation \r{D22-1}.

\section*{Conclusion}

This paper investigates quantum loop algebra $U_q(\Dn{n})$.
We discuss the connection between two realizations of this algebra: one in terms of the fundamental $\LL$-operators and Jimbo $\RR$-matrix \cite{Jimbo1986quantum} and the second in terms of the currents \cite{Drinfeld1988NewRealisation}.
It was shown that in order to apply a standard approach to the description of this algebra in terms of the Gaussian coordinates of the fundamental $\LL$-operators, one has to slightly modify Jimbo's $\RR$-matrix \r{RJmat}.

One can use these results to describe Bethe vectors for the quantum integrable models associated with quantum loop algebra $U_q(\Dn{n})$ as it was done in \cite{KyotoPaper, HutsalyukLiashykPakuliakRagoucySlavnov2017current}.
Extending the results obtained in the paper to exceptional quantum loop algebras is an interesting open problem.

\section*{Acknowledgement}
S.P. acknowledges the support of the PAUSE Programme  and hospitality at LAPTh where
this work was finalized.
The research of A.L. was supported by Beijing Natural Science Foundation (IS24006).
This work was partially done when A.L. visited the Max Plank Institut f\"ur Mathematik in Bonn.
He thanks the MPIM for the hospitality and stimulating scientific atmosphere.
A.L. is also grateful to the CNRS PHYSIQUE for support during his visit to Annecy in the course of this investigation.

\appendix

\section{Jimbo's $\RR$-matrix}\label{ApA}

In \cite{Jimbo1986quantum} $\RR$-matrix
$\RR^{\rm J}(u,v)\in{\textrm{End}}(\CC^{N}\ot\CC^{N})$
 associated with the vector
representation of $\mathfrak{o}_{2n}$ was presented in the form
\begin{equation}\label{RJmat}
\begin{split}
&\RR^{\rm J}(u,v)\ =\
f(u^2,v^2) \sum_{i \ne n, n+1} \E_{ii}\ot \E_{ii} +
\sum_{i \ne j,j' \atop i \text{ or } j \text{ }\ne\text{ } n,n+1}
\E_{ii}\ot \E_{jj}\\
&+\sum_{i<j \atop i,j \ne n, n+1}
\left(\gle(u^2,v^2) \E_{ij}\ot \E_{ji}+ \gri(u^2,v^2) \E_{ji}\ot \E_{ij}\right)
+\sum_{i,j \ne n,n+1} a_{ij}(u,v) \ \E_{i'j'}\ot \E_{ij}\\
&+\frac{1}{2}\sum_{1 \le i \le n-1 \atop j = n, n+1}
\Big(\gle(u,v) \left( \E_{ij}\ot \E_{ji} + \E_{j i'}\ot \E_{i' j} \right)
+ \gri(u,v) \left( \E_{ji}\ot \E_{ij} + \E_{i'j}\ot \E_{ji'} \right) \\
&\quad+
\gle(-u,v) \left( \E_{ij}\ot \E_{j'i} + \E_{ji'}\ot \E_{i'j'}\right)  +
\gri(-u,v) \left( \E_{ji}\ot \E_{ij'} + \E_{i'j}\ot \E_{j'i'}\right)
\Big) \\
&+\frac{1}{2}\sum_{ i \ne n,n+1 \atop j = n, n+1}
\Big( b^{+}_i(u,v) \left( \E_{ij}\ot \E_{i'j'} + \E_{j'i'}\ot \E_{ji} \right) +
 b^{-}_i(u,v) \left( \E_{ij}\ot \E_{i'j} + \E_{ji'}\ot \E_{ji} \right)
\Big) \\
&+ \sum_{i=n,n+1} \Big(
 c^{+}(u,v)\ \E_{ii}\ot \E_{i'i'} +
 c^{-}(u,v)\ \E_{ii}\ot \E_{ii} +\\
 &\qquad\qquad\qquad\qquad\qquad +d^{+}(u,v)\ \E_{i'i}\ot \E_{ii'} +
 d^{-}(u,v)\ \E_{ii'}\ot \E_{ii'}
 \Big),
 \end{split}
\end{equation}
where
\begin{subequations}\label{funct}
\begin{equation}\label{functa}
a_{ij}(u,v)=\begin{cases}
f(v^2\xi^2 , u^2 ),\quad &i=j,\\
q^{\bar\imath-\bar\jmath}
\gle(v^2\xi^2,u^2),\quad &i<j,\\
q^{\bar\imath-\bar\jmath}
\gri(v^2\xi^2,u^2),\quad &i>j,
\end{cases}
\end{equation}
\begin{equation}\label{functb}
    b^{\pm}_i(u,v) =
\begin{cases}
 -q^{\bar n-\bar\imath } \gle(\pm u,v\xi), \quad i<n,\\
- q^{\bar n-\bar\imath} \gri(\pm u,v\xi), \quad i>n+1,
\end{cases}
\end{equation}
\begin{equation}\label{functc}
    c^{\pm}(u,v) = 1 + \frac12\frac{1+\xi}{q - q^{-1}} \gri(\mp u,v) \gle(\pm u,v\xi),
    \end{equation}
    \begin{equation}\label{functd}
    d^{\pm}(u,v) = \frac12\frac{1-\xi}{q - q^{-1}} \gri(\pm u,v) \gle(\pm u,v\xi).
\end{equation}
\end{subequations}

In order to identify the expression for $\RR$-matrix given by \r{RJmat} with
the one presented in \cite{Jimbo1986quantum}, one has to do the following. First,
one should shift $n\to n-1$ in (3.7) of \cite{Jimbo1986quantum}. Then one has to
impose following identification
\begin{equation*}
x\to \frac{v}{u},\qquad k\to q,\qquad \xi\to \xi^{-1}.
\end{equation*}
The map $\bar\alpha$ used in \cite{Jimbo1986quantum} after shifting $n\to n-1$ is related to the map $\bar\imath$ \r{map} by the relation
\begin{equation*}
\bar\alpha=-\bar\imath+n+\frac12.
\end{equation*}
One can check that using this identification of the parameters expression (3.7)
from \cite{Jimbo1986quantum} coincides identically with \r{RJmat}.

Using explicit form of the operator $\Ut$ \r{Ut} and
the fact that $n'=n+1$ and $(n+1)'=n$ one can obtain
\begin{equation}\label{A1}
\Ut\cdot\E_{ij}\cdot \Ut=\E_{ij}\ \delta_{i,j\not=n,n+1}+
\E_{ij'}\ \delta_{i\not=n,n+1\atop j=n,n+1}+
\E_{i'j}\ \delta_{i=n,n+1\atop j\not=n,n+1}+
\E_{i'j'}\ \delta_{i=n,n+1\atop j=n,n+1}
\end{equation}
and write
\begin{equation}\label{A2}
1= \delta_{i,j\not=n,n+1}+
\delta_{i\not=n,n+1\atop j=n,n+1}+
\delta_{i=n,n+1\atop j\not=n,n+1}+
\delta_{i=n,n+1\atop j=n,n+1}
\end{equation}
for any $i,j=1,\ldots,2n$.

One can observe that
\begin{equation}\label{A3}
\Id\ot\Id=\sum_{i\not=n,n+1}\E_{ii}\ot\E_{ii}+\sum_{i\not=n,n+1}\E_{ii}\ot\E_{i'i'}+
\sum_{i,j=n,n+1}\E_{ii}\ot\E_{jj}+
\sum_{i\not=j,j'\atop i\,{\rm or}\,j\not=n,n+1}\E_{ii}\ot\E_{jj}\,.
\end{equation}

Let us denote by $\PPd(u,v)$ and $\QQd^{\rm J}(u,v)$ the matrices
\begin{equation}\label{PP14}
\PPd(u,v)=\frac12\Big(\Per(u,v)+(\Ut\ot\Id)\cdot \Per(-u,v)\cdot (\Ut\ot\Id)\Big)\,,
\end{equation}
\begin{equation}\label{QQd}
\QQd^{\rm J}(u,v)=\frac12\Big(\Qer^{\rm J}(u,v)+(\Ut\ot\Id)\cdot \Qer^{\rm J}(-u,v)\cdot (\Ut\ot\Id)\Big)\,.
\end{equation}

Then $\RR$-matrix $\RR(u,v)$ \r{D2Rmat} takes the form
\begin{equation}\label{D2Ralt}
 \RR^{\rm J}(u,v)=\Id\ot\Id+\PPd(u,v)+\QQd^{\rm J}(u,v)\,.
\end{equation}

Using \r{A1} and \r{A2} operator $\PPd(u,v)$ can be written in the form
\begin{equation}\label{A4}
\begin{aligned}
\PPd(u,v)=\frac12 \sum_{1\leq {i,j}\leq 2n}&\left[\Big(\pf_{ij}(u,v)+\pf_{ij}(-u,v)\Big)\ \E_{ji}\ot\E_{ij}\
\delta_{i,j\not=n,n+1}+\right.\\
+ &\Big(\pf_{ij}(u,v)\ \E_{ji}+\pf_{ij}(-u,v)\ \E_{j'i}  \Big)\ot \E_{ij}\
\delta_{ i\not=n,n+1\atop j=n,n+1}+\\
+ &\Big(\pf_{ij}(u,v)\ \E_{ji}+\pf_{ij}(-u,v)\ \E_{ji'}  \Big)\ot \E_{ij}\
\delta_{ i=n,n+1\atop j\not=n,n+1}+\\
+ &\left.\Big(\pf_{ij}(u,v)\ \E_{ji}+\pf_{ij}(-u,v)\ \E_{j'i'}  \Big)\ot \E_{ij}\
\delta_{i,j=n,n+1}\right]\,.
\end{aligned}
\end{equation}
Using trivial identities $f(u,v)+f(-u,v)=2f(u^2,v^2)$ and
\begin{equation}\label{sum-eq}
 \gle(u,v)+\gle(-u,v)=2\gle(u^2,v^2),\qquad  \gri(u,v)+\gri(-u,v)=2\gri(u^2,v^2)
\end{equation}
the first line in \r{A4} yields
\begin{equation*}
(f(u^2,v^2)-1)\sum_{i\not=n,n+1}\E_{ii}\ot\E_{ii}+
\sum_{i<j\atop i,j\not=n,n+1}(\gle(u^2,v^2)\E_{ij}\ot\E_{ji}+
\gri(u^2,v^2)\E_{ji}\ot\E_{ij})
\end{equation*}
which together with first and fourth terms in \r{A3} yields first three terms in \r{RJmat}.

The second and third lines of \r{A4} result in the third and fourth lines of \r{RJmat}.
To obtain in \r{RJmat} the terms proportional to the functions $a_{ij}(u,v)$,
$b^\pm_i(u,v)$, $c^\pm(u,v)$ and $d^\pm(u,v)$ given by equalities \r{funct}
one has to present the operators $\QQd(u,v)$ in the form
\begin{equation}\label{A5}
\begin{split}
\QQd^{\rm J}(u,v)=\frac12
 \sum_{1\leq {i,j}\leq 2n}&\left[\Big(\qf_{ij}(u,v)\ \E_{i'j'}+\qf_{ij}(-u,v)\ \E_{i'j'}\Big)\ \ot\E_{ij}\
\delta_{i,j\not=n,n+1}+\right.\\
+ &\Big(\qf_{ij}(u,v)\ \E_{i'j'}+\qf_{ij}(-u,v)\ \E_{i'j}  \Big)\ot \E_{ij}\
\delta_{i\not=n,n+1\atop j=n,n+1}+\\
+ &\Big(\qf_{ij}(u,v)\ \E_{i'j'}+\qf_{ij}(-u,v)\ \E_{ij'}  \Big)\ot \E_{ij}\
\delta_{i=n,n+1\atop j\not=n,n+1}+\\
+ &\left.\Big(\qf_{ij}(u,v)\ \E_{i'j'}+\qf_{ij}(-u,v)\ \E_{ij}  \Big)\ot \E_{ij}\
\delta_{i,j=n,n+1}\right]\,.
\end{split}
\end{equation}
Now the first line in \r{A5} together with the second terms in \r{A3} yields the definition
of the function $a_{ij}(u,v)$ given by \r{functa}. Second and third lines in \r{A5} produces
the terms in \r{RJmat} proportional to the functions $b^\pm_i(u,v)$ given by the equality
\r{functb}. Finally the forth lines in \r{A4} and \r{A5} together with third term from
\r{A3} yields the terms
proportional to the functions $c^\pm(u,v)$ \r{functc} and $d^\pm(u,v)$ \r{functd}
in \r{RJmat}.

\subsection{``Diagonal'' Gaussian coordinates for Jimbo's $\RR$-matrix}
\label{sec:GD Jimbo R-matrix}

One can define also quantum loop algebra associated with the
Jimbo's $\RR$-matrix $\RR^{\rm J}(u,v)$ \r{D2Rmat} and introduce
similar Gaussian coordinates. Matrix entries $\FF^\pm_{j,n}(u)$,
$\FF^\pm_{j,n+1}(u)$, $\FF^\pm_{n,i}(u)$, $\FF^\pm_{n+1,i}(u)$ and
$\EE^\pm_{n,j}(u)$,
$\EE^\pm_{n+1,j}(u)$, $\EE^\pm_{i,n}(u)$, $\EE^\pm_{i,n+1}(u)$
will be replaced by the linear combinations of the same Gaussian coordinates
in upper and lower  triangular
matrices $\mathbf{F}^\pm(u)$, $\mathbf{E}^\pm(u)$.
The matrix $\mathbf{K}^\pm(u)$ becomes
 'almost' diagonal
\begin{equation}\label{ndK}
\tilde{\mathbf{K}}^\pm(u)= \sum_{1\leq i\leq 2n} \E_{ii}\ k^\pm_{i}(u)
+\sum_{i=n,n+1}\E_{i'i}\ k^\pm_{i,i'}(u)\,.
\end{equation}

Due to the relation \r{tRalt}, it is obvious
that the relation  between $\LL$-operators $\tilde{\LL}^\pm(u)$ for
Jimbo's $\RR$-matrix $\RR^{\rm J}(u,v)$ and with those given by \r{Lop} is
\begin{equation}\label{LtL}
\tilde{\LL}^\pm(u)=\TT^{-1}\ \LL^\pm(u)\ \TT
\end{equation}
which results in the relations
\begin{equation}\label{ktk}
\begin{split}
  k^\pm_{n,n}(u)=k^\pm_{n+1,n+1}(-u)&= \frac{(\xi^{1/2} + 1)^{2}}{4\, \xi^{1/2}} k^\pm_n(u) - \frac{(\xi^{1/2}-1)^{2}}{4\, \xi^{1/2}} k^\pm_n(-u)\,,\\
k^\pm_{n,n+1}(u)=-k^\pm_{n+1,n}(u)&=
\frac{1 - \xi}{4\, \xi^{1/2}}\Big(k^\pm_n(u)-k^\pm_n(-u)\Big)\,.
\end{split}
\end{equation}

\section{Crossing symmetries and pole structure}
\label{sec: Crossing symmetries and pole structure}

By definition of the matrices $\Per(u,v)$ \r{PPuv} and $\Qer^{\rm J}(u,v)$
\r{QQuv} and property
\begin{equation*}
\pf_{ij}(u,v)=\pf_{j'i'}(u,v),\quad \qf_{ij}(u,v)=\qf_{j'i'}(u,v)
\end{equation*}
one gets
\begin{equation*}
\Per(u,v)^{{\rm t}_1{\rm t}_2}=\Per(u,v),\quad
\Qer(u,v)^{{\rm t}_1{\rm t}_2}=\Qer(u,v)\,.
\end{equation*}
It implies  that
\begin{equation}\label{reflecJ}
\RR^{\rm J}_{12}(u,v)^{{\rm t}_1{\rm t}_2}=\RR^{\rm J}_{12}(u,v).
\end{equation}
Using antimorphism property of the transposition, one gets from \r{reflecJ}
the property \r{reflec} for the $\RR$-matrix $\RR(u,v)$.

Due to the commutativity $[\Ut,\TT]=[\Der,\TT]=0$,
 $\RR$-matrix $\RR^{\rm J}(u,v)$ has the same properties  \r{comm1}
 and \r{comm2} as  $\RR$-matrix $\RR(u,v)$.

 Using relation \r{con1} one can find that  $\RR$-matrix $\RR^{\rm J}(u,v)$
 satisfies an equality
\begin{equation}\label{crre1}
\Der_2\ \Per_{12}\ \RR^{\rm J}_{12}(v\xi,u)^{{\rm t}_1}\ \Per_{12}\ \Der_2^{-1}=
\RR^{\rm J}_{12}(u,v)
\end{equation}
 or
\begin{equation}\label{crre2}
\Der_1\  \RR^{\rm J}_{12}(v\xi,u)^{{\rm t}_1}\  \Der_1^{-1}=
\RR^{\rm J}_{21}(u,v)\,.
\end{equation}
Conjugating the latter equality by $\TT_1\TT_2$ one gets \r{cross1}
or using unitarity \r{unitar} in another form
 \begin{equation}\label{crre3}
\Der_1\ \TT_1^2\ \RR_{12}(v\xi,u)^{{\rm t}_1}\ \TT_1^{-2}\  \Der_1^{-1}
\RR_{12}(v,u)=f(u^2,v^2)f(v^2,u^2)\Id\ot\Id\,.
\end{equation}
Considering equality \r{crre3} at $v\to v\xi$ and conjugating left hand side
of this relation by $\Der_1\ \TT_1^{-2}$ one gets
\begin{equation}\label{crre4}
\begin{split}
&\Big(\Der^2_1\  \RR_{12}(v\xi^2,u)^{{\rm t}_1}\   \Der_1^{-2} \Big)\cdot
\Big(\Der_1\ \TT_1^{-2}\ \RR_{12}(v\xi,u)\ \TT_1^{2}\ \Der_1^{-1}\Big)=\\
&\qquad\qquad= f(u^2,v^2\xi^2)f(v^2\xi^2,u^2)\Id\ot\Id\,.
\end{split}
\end{equation}
Taking transposition of equality \r{crre2} in the first tensor space and substituting it
into \r{crre4} one proves \r{cross2}.

Now one can calculate the residues of the
$\RR$-matrix $\RR^{\rm J}(u,v)$ at the poles $u=\pm v$ and $u=\pm v\xi$.
It is clear that only term $\Per(u,v)$ will contribute to the residue at $u=v$ and
\begin{equation}\label{scre5}
\left.\frac{2 (u-v)}{u (q-q^{-1})}\ \RR^{\rm J}_{12}(u,v)\right|_{u=v}=\Per_{12}\,.
\end{equation}
Then residue at the pole $u=v\xi$ can be calculated using \r{scre5} as follows
\r{crre1} to obtain
\begin{equation}\label{scre6}
\left.\frac{2 (v\xi-u)}{u (q-q^{-1})}\ \RR^{\rm J}_{12}(v\xi,u)\right|_{u=v}=
\Der_1^{-1}\ \Per^{{\rm t}_1}_{12}\ \Der_1\,.
\end{equation}
To obtain right hand side in \r{scre6}, one has to use an equality
$\Der_2\ \Per_{12}^{{\rm t}_1}=\Der^{{\rm t}_1}_1\ \Per_{12}^{{\rm t}_1}=
\Der^{-1}_1\ \Per_{12}^{{\rm t}_1}$ since $\Der^{\rm t}=\Der^{-1}$.
Conjugating both sides of \r{scre6} by $\TT_1\ \TT_2$ and using
$\TT_2\ \Per_{12}^{{\rm t}_1}=\TT^{{\rm t}_1}_1\ \Per_{12}^{{\rm t}_1}=
\TT_1\ \Per_{12}^{{\rm t}_1}$ where $\TT^{\rm t}=\TT$ one obtained third
residue in \r{poles}. Other residues in \r{poles} can be considered
analogously.

To calculate residue in \r{cen2} one can express $\Big(\RR_{12}(u,v)^{{\rm t}_1}\Big)^{-1}$
from \r{cross2}
\begin{equation*}
f(u^2,v^2\xi^2)f(v^2\xi^2,u^2)\Big(\RR_{12}(u,v)^{{\rm t}_1}\Big)^{-1}=
\Der_2^2\ \RR_{21}(v\xi^2,u)^{{\rm t}_2}\ \Der_2^{-2}
\end{equation*}
and then use \r{scre5}.

\subsection{Another presentation of the matrix $\Qer(u,v|\xi)$}

In this appendix, we will prove following proposition. 

\begin{prop}\label{Qprop}
In representation \eqref{tRalt} the matrix $\Qer(u,v|\xi)$ \r{Qeven} can be replaced by the matrix $\tilde\Qer(u,v|\xi)$, which is explicitly even in the parameter $\xi$ 
\begin{equation}\label{Qev}
\tilde\Qer(u,v|\xi)=\tilde\Qer(u,v|-\xi)
\end{equation}
and has the form
\begin{equation}\label{Qzero}
\begin{split}
&\tilde\Qer(u,v|\xi)=\frac12\sum_{1\leq i,j\leq N\atop i,j\not=n,n+1}
(\qf_{ij}(u,v|\xi)+\qf_{ij}(u,v|-\xi))\E_{i'j'}\ot\E_{ij}+\\
&\quad+\frac{1}{2q^{1/2}}
\sum_{1\leq j\leq N\atop {i=n,n+1\atop j\not=n,n+1}}
\Big((1+\xi^{-1})\qf_{ij}(u,v|\xi)+(1-\xi^{-1})\qf_{ij}(u,v|-\xi)\Big)\E_{i'j'}\ot\E_{ij}+\\
&\quad+\frac{1}{2q^{-1/2}}
\sum_{1\leq i\leq N\atop {i\not=n,n+1\atop j=n,n+1}}
\Big((1+\xi)\qf_{ij}(u,v|\xi)+(1-\xi)\qf_{ij}(u,v|-\xi)\Big)\E_{i'j'}\ot\E_{ij}+\\
&\quad+\frac14\sum_{i,j=n,n+1}\Big(
\QQff(u,v|\xi) \E_{i'j'}\ot\E_{ij}+
\QQdd(u,v|\xi)(\E_{ij}-\E_{i'j'})\ot(\E_{ij'}-\E_{i'j})\Big)-\\
&-\frac14\Big(\alpha_q(\E_{n,n}+\E_{n+1,n+1})\ot (\E_{n,n}+\E_{n+1,n+1})+
\gamma_q(\E_{n,n+1}+\E_{n+1,n})\ot (\E_{n+1,n}-\E_{n,n+1})\Big)\,,
\end{split}
\end{equation}
where 
\begin{equation}\label{QQff}
\QQff(u,v|\xi)=\gamma_q\frac{(v\xi^2+ u)(v+ u)}{v^2\xi^2-u^2}, \quad
\QQdd(u,v|\xi)=\gamma_q\frac{(\xi^2-1)\, v\, u}
{v^2\xi^2-u^2}\,.
\end{equation}
\end{prop}

We will calculate the matrix $\Qer(u,v|\xi)$ by the explicit conjugation of 
$\RR^{\rm J}(u,v)$ using presentation of the matrices 
$\Per(u,v)$ and  $\Qer^{\rm J}(u,v|\xi)$ given by equalities \r{A4} and \r{A5}.

The inverse matrix $\TT^{-1}$ is 
 \begin{equation}\label{Tm2}
\TT^{-1}_{ij}=\delta_{ij}\delta_{i,j\not=n,n+1}+\big(a_+'\delta_{ij}+a_-'\delta_{ij'})
\delta_{i,j=n,n+1}\,,
\quad a_\pm'=\frac{1\pm\xi^{1/2}}{2q^{-1/4}}\,.
\end{equation}
Using \r{Tm1} and \r{Tm2} one gets 
\begin{equation}\label{Tm3}
\begin{split}
\TT\ \E_{ij}\ \TT^{-1}&=\E_{ij}\delta_{i,j\not=n,n+1}+
(a_+\E_{ij}+a_-\E_{i'j})\delta_{i=n,n+1\atop j\not=n,n+1}+
(a_+'\E_{ij}+a_-'\E_{ij'})\delta_{i\not=n,n+1\atop j=n,n+1}+\\
&+\Big(a_+a_+'\E_{ij}+a_+'a_-\E_{i'j}+a_+a_-'\E_{ij'}+a_-a_-'\E_{i'j'}\Big)\delta_{i,j=n,n+1}\,.
\end{split}
\end{equation}

Using the first three terms in the right hand side of \r{Tm3}, one can observe that 
first three lines of \r{A4} do not change after conjugation of 
$\RR^{\rm J}(u,v)$ in \r{tRalt}.
Analogously, one can verify that conjugation of the first three lines in \r{A5} 
produces first three lines in \r{Qzero}. 

For $i,j=n,n+1$ functions \r{p-fun} and \r{q-fun} has the presentation
which separates dependence of these functions on spectral  parameters and 
indices $i,j$
\begin{equation}\label{pq-fun-alt}
\begin{split}
\pf_{ij}(u,v)&=\frac{\Gap(u,v)}{2}+\frac{\alpha_q\delta_{ij}}{2}+
\frac{\epsilon_j\,\gamma_q\,\delta_{ij'}}{2}\,,\\
\qf_{ij}(u,v|\xi)&=\frac{\Gap(v\xi,u)}{2}-\frac{\alpha_q\delta_{ij}}{2}-
\frac{\epsilon_j\,\gamma_q\,\delta_{ij'}}{2}\,,
\end{split}
\end{equation} 
where 
\begin{equation}\label{pi-fun}
\Gap(u,v)=\gle(u,v)+\gri(u,v)=\gamma_q\, \frac{u+v}{u-v},\quad  
\epsilon_j = (-1)^{j - n}\,.
\end{equation}

To finish the calculation of the conjugation $\TT_1\TT_2\,\RR^{\rm J}(u,v)\,\TT_2^{-1}
\TT_1^{-1}$ one has to  sum over $i,j=n,n+1$ 
the terms 
\begin{equation}\label{Tm6}
\begin{split}
X(u,v)=
\frac12\sum_{i,j=n,n+1}
\Big(&\pf_{ij}(u,v)\TT\E_{ji}\TT^{-1}\ot \TT\E_{ij}\TT^{-1}+
\pf_{ij}(-u,v)\TT\E_{j'i'}\TT^{-1}\ot \TT\E_{ij}\TT^{-1}+\\
+&\qf_{ij}(u,v|\xi)\TT\E_{i'j'}\TT^{-1}\ot \TT\E_{ij}\TT^{-1}+
\qf_{ij}(-u,v|\xi)\TT\E_{ij}\TT^{-1}\ot \TT\E_{ij}\TT^{-1}
\Big),
\end{split}
\end{equation} 
where we can replace the functions 
$\pf_{ij}(u,v)$ and $\qf_{ij}(u,v|\xi)$  using formulas 
\r{pq-fun-alt}. One can verify that terms in \r{Tm6} which do not depend on the 
spectral parameters will cancel each other after summation over $i,j=n,n+1$ and 
matrix $X(u,v)$ is equal to 
\begin{equation}\label{Tm61}
\begin{split}
X(u,v)=
\frac14\sum_{i,j=n,n+1}
\Big(&\Gap(u,v)\E_{ji}\ot \E_{ij}+
\Gap(-u,v)\E_{j'i'}\ot \E_{ij}+\\
+&\Gap(v\xi,u)\TT\E_{i'j'}\TT^{-1}\ot \TT\E_{ij}\TT^{-1}+
\Gap(v\xi,-u)\TT\E_{ij}\TT^{-1}\ot \TT\E_{ij}\TT^{-1}
\Big)\,,
\end{split}
\end{equation} 
where in the first line of \r{Tm61} we used simple identities 
\begin{equation}
\begin{split}
&\sum_{i,j=n,n+1}\Big(\TT\E_{ji}\TT^{-1}\ot \TT\E_{ij}\TT^{-1}-\E_{ji}\ot\E_{ij}\Big)=0\,,\\
&\sum_{i,j=n,n+1}\Big(\TT\E_{j'i'}\TT^{-1}\ot \TT\E_{ij}\TT^{-1}-\E_{j'i'}\ot\E_{ij}\Big)=0\,.
\end{split}
\end{equation}
The first identity is a consequence of  $\TT_1\TT_2\Per_{12}=\Per_{12}\TT_1\TT_2$
and the second identity follows from the first one due to the commutativity 
$[\TT,\Ut]=0$. 

Using \r{pq-fun-alt} once again we can rewrite the matrix $X(u,v)$ in the form 
\begin{equation}\label{Tm62}
\begin{split}
X(u,v)&=\frac12\sum_{i,j=n,n+1}\Big(\pf_{ij}(u,v)\E_{ji}+\pf_{ij}(-u,v)\E_{j'i'}\Big)\ot\E_{ij}+\\
&\quad+\frac14\sum_{i,j=n,n+1}\Big(\Gap(v\xi,u)\TT\E_{i'j'}\TT^{-1}\ot \TT\E_{ij}\TT^{-1}+
\Gap(v\xi,-u)\TT\E_{ij}\TT^{-1}\ot \TT\E_{ij}\TT^{-1}
\Big)-\\
&\quad-\frac14\sum_{i,j=n,n+1}(\alpha_q\,\delta_{ij}+\epsilon_j\,\gamma_q\,
\delta_{ij'})(\E_{ji}+\E_{j'i'})\ot\E_{ij}\,.
\end{split}
\end{equation}

The terms from the first line in \r{Tm62} will restore the matrix 
\begin{equation*}
\PPd(u,v)=\frac12\Big(\Per(u,v)+(\Ut\ot\Id)\cdot\Per(-u,v)\cdot(\Ut\ot\Id)\Big)
\end{equation*} 
in \r{tRalt}.  

Let us calculate 
\begin{equation}\label{TM63}
\sum_{i,j=n,n+1}\TT\E_{i'j'}\TT^{-1}\ot \TT\E_{ij}\TT^{-1}
\quad\mbox{and}
\sum_{i,j=n,n+1}\TT\E_{ij}\TT^{-1}\ot \TT\E_{ij}\TT^{-1}
\end{equation}
from the second line of \r{Tm62} using the fourth term in the 
conjugation \r{Tm3} of the matrix $\E_{ij}$.
Due to the commutativity $[\TT,\Ut]=0$ it is sufficient to calculate 
first formula in \r{TM63}. The second one can be obtained by the conjugation 
with the matrix $\Ut$ in the first tensor component.

Due to the trivial equalities 
\begin{equation*}
a_\pm a_\pm'=\frac{1\pm\xi_+}{2},\quad 
a_+'a_-=-a_+a_-'=\frac{\xi_-}{2},\quad \xi_\pm=\frac{\xi^{1/2}\pm\xi^{-1/2}}{2}\,.
\end{equation*}
one can note that for $i,j=n,n+1$
\begin{equation}\label{Tm4}
\TT\,\E_{ij}\,\TT^{-1}=
\frac12\Big(\adm_{ij}+\xi_+\bdm_{ij}+\xi_-\bdm_{i'j}\Big)\,,
\end{equation}
where 
\begin{equation}\label{Tm5}
\adm_{ij}=\E_{ij}+\E_{i'j'}=\adm_{i'j'}\,,\quad
\bdm_{ij}=\E_{ij}-\E_{i'j'}=-\bdm_{i'j'}\,,\quad i,j=n,n+1\,.
\end{equation}
Using \r{Tm4} and the properties \r{Tm4} one can calculate 
\begin{equation}\label{tm1}
\begin{split}
&\sum_{i,j=n,n+1}\TT\E_{i'j'}\TT^{-1}\ot \TT\E_{ij}\TT^{-1}
=\\
&\qquad=\sum_{i,j=n,n+1}\sk{\frac{(1+\xi)^2}{16\xi}\E_{i'j'}\ot \E_{ij}-
\frac{(1-\xi)^2}{16\xi}\E_{ij}\ot \E_{ij}+
\frac{\xi^2-1}{16\xi}\bdm_{ij}\ot\bdm_{ij'}}
\end{split}
\end{equation}
and 
\begin{equation}\label{tm2}
\begin{split}
&\sum_{i,j=n,n+1}\TT\E_{ij}\TT^{-1}\ot \TT\E_{ij}\TT^{-1}
=\\
&\qquad=\sum_{i,j=n,n+1}\sk{\frac{(1+\xi)^2}{16\xi}\E_{ij}\ot \E_{ij}-
\frac{(1-\xi)^2}{16\xi}\E_{i'j'}\ot \E_{ij}-
\frac{\xi^2-1}{16\xi}\bdm_{ij}\ot\bdm_{ij'}}
\end{split}
\end{equation}
Substituting \r{tm1} and \r{tm2} into second line of \r{Tm62} we prove the 
statement of proposition~\ref{Qprop} due to the identities for the rational 
functions
\begin{equation*}
\begin{split}
\frac{1}{4\xi}\Big((1+\xi)^2\Gap(v\xi,u)-(1-\xi)^2\Gap(v\xi,-u)\Big)&=
\gamma_q\frac{(v\xi^2+ u)(v+ u)}{v^2\xi^2-u^2}, \\
\frac{1}{4}\Big(\Gap(v\xi,u)-\Gap(v\xi,-u)\Big)&=
\gamma_q\frac{ v\xi\, u}
{v^2\xi^2-u^2}\,.
\end{split}
\end{equation*}
Last line in \r{Tm62} is equal to the constant term of the matrix 
$\tilde\Qer(u,v|\xi)$ given by the last line in \r{Qzero}.
\qed

\section{Proof of lemma~\ref{5.3}}\label{ApD}

To prove this lemma, one needs one more lemma, which claims the identities
for the rational functions.
\begin{lemma}\label{iden}
There are identities for the rational functions with $\ell, j > 1$
\begin{equation}\label{em19}
\begin{split}
&\big(\qf_{\ell j}(u,v|\xi)+
\qf_{\ell j}(-u,v|\xi)\big)-\\
&\qquad- \frac{q}{2}
\big(\qf_{1 j}(u,v|\xi)+
\qf_{1 j}(-u,v|\xi)\big)
\big(\qf_{\ell 1}(u,qv|\xi)+
\qf_{\ell 1}(-u,qv|\xi)\big)=\\
&\qquad=\big(\qf_{\ell j}(u,v|q\xi)+
\qf_{\ell j}(-u,v|q\xi)\big)\,,\quad\mbox{for any}\quad\ell,j,
\end{split}
\end{equation}
\begin{equation}\label{em20}
\begin{split}
&\big(\qf_{\ell j}(u,v|\xi)-
\qf_{\ell j}(-u,v|\xi)\big)-\\
&\qquad- \frac{1}{2}
\big(\qf_{1 j}(u,v|\xi)+
\qf_{1 j}(-u,v|\xi)\big)
\big(\qf_{\ell 1}(u,qv|\xi)-
\qf_{\ell 1}(-u,qv|\xi)\big)=\\
&\qquad=q^{-1}\big(\qf_{\ell j}(u,v|q\xi)-
\qf_{\ell j}(-u,v|q\xi)\big),\quad\mbox{for}\quad\ell\not=j,
\end{split}
\end{equation}
\begin{equation}\label{em18}
\begin{split}
&\big(\qf_{\ell j}(u,v|\xi)-\qf_{\ell j}(-u,v|\xi)\big)-\\
&\qquad- \frac{q}{2}
\big(\qf_{1 j}(u,v|\xi)-\qf_{1 j}(-u,v|\xi)\big)
\big(\qf_{\ell 1}(u,qv|\xi)+\qf_{\ell 1}(-u,qv|\xi)\big)=\\
&\qquad=q\big(\qf_{\ell j}(u,v|q\xi)-
\qf_{\ell j}(-u,v|q\xi)\big),\quad\mbox{for}\quad\ell\not=j.
\end{split}
\end{equation}

\end{lemma}

To prove \r{em19} one has to consider separately three cases $\ell<j$, $\ell>j$ and
$\ell=j$. In all three cases, right hand side of \r{em19} will be proportional to
$q^{\bar\ell-\bar\jmath}$. Then according to \r{q-fun} for $\ell<j$
\begin{equation*}
\begin{split}
&\big(\gle(v\xi,u)+\gle(v\xi,-u)\big)-\frac{q}{2}
\big(\gle(v\xi,u)+\gle(v\xi,-u)\big)
\big(\gri(qv\xi,u)+\gri(qv\xi,-u)\big)=\\
&\quad=\frac{2(q-q^{-1})v^2\xi^2}{v^2\xi^2-u^2}\Big(1-
\frac{(q-q^{-1})u^2}{q^2v^2\xi^2-u^2}\Big)=\frac{2(q-q^{-1})q^2v^2\xi^2}{q^2v^2\xi^2-u^2}=\\
&\quad=\Big(\frac{(q-q^{-1})qv\xi}{qv\xi-u}+\frac{(q-q^{-1})qv\xi}{qv\xi-(-u)}\Big)=
\gle(qv\xi,u)+\gle(qv\xi,-u)\,.
\end{split}
\end{equation*}
Analogously for $\ell>j$
\begin{equation*}
\begin{split}
&\big(\gri(v\xi,u)+\gri(v\xi,-u)\big)-\frac{q}{2}
\big(\gle(v\xi,u)+\gle(v\xi,-u)\big)
\big(\gri(qv\xi,u)+\gri(qv\xi,-u)\big)=\\
&\quad=\frac{2(q-q^{-1})u^2}{v^2\xi^2-u^2}\Big(1-
\frac{(q-q^{-1})v^2\xi^2}{q^2v^2\xi^2-u^2}\Big)=\frac{2(q-q^{-1})u^2}{q^2v^2\xi^2-u^2}=\\
&\quad=\Big(\frac{(q-q^{-1})u}{qv\xi-u}+\frac{(q-q^{-1})(-u)}{qv\xi-(-u)}\Big)=
\gri(qv\xi,u)+\gri(qv\xi,-u)\,.
\end{split}
\end{equation*}
Finally for $\ell=j \ne n, n+1$
\begin{equation*}
\begin{split}
&\big(f(v\xi,u)+f(v\xi,-u)\big)-\frac{q}{2}
\big(\gle(v\xi,u)+\gle(v\xi,-u)\big)
\big(\gri(qv\xi,u)+\gri(qv\xi,-u)\big)=\\
&\quad= \frac{2}{v^2\xi^2-u^2}\Big(qv^2\xi^2-q^{-1}u^2-
q\frac{(q-q^{-1})^2v^2\xi^2u^2}{q^2v^2\xi^2-u^2}\Big)=
2\frac{q^3v^2\xi^2-q^{-1}u^2}{q^2v^2\xi^2-u^2}=\\
&\quad=\Big(\frac{q^2v\xi-q^{-1}u}{qv\xi-u}+\frac{(q^2v\xi-(q^{-1}(-u)}{qv\xi-(-u)}\Big)=
f(qv\xi,u)+f(qv\xi,-u)\,.
\end{split}
\end{equation*}
For the case $\ell = j = n, n+1$ proving is the same up to addition by the same constant in both sides of the relation.

To prove \r{em20} we again observe that it is proportional to $q^{\bar\ell-\bar\jmath}$.
Then one can note that
\begin{equation}\label{g-g}
\gle(v\xi,u)-\gle(v\xi,-u)=\gri(v\xi,u)-\gri(v\xi,-u)=\frac{2(q-q^{-1})v\xi u}{v^2\xi^2-u^2}
\end{equation}
and for both cases $\ell<j$ or $\ell>j$ the right hand side of \r{em20} is equal to
\begin{equation*}
\begin{split}
&\frac{2(q-q^{-1})v\xi u}{v^2\xi^2-u^2}-
\frac{2(q-q^{-1})v^2\xi^2}{v^2\xi^2-u^2}\
\frac{(q-q^{-1})q v\xi u}{q^2v^2\xi^2-u^2}=
\frac{2(q-q^{-1})v\xi u}{q^2v^2\xi^2-u^2}=\\
&\quad=\frac{(q-q^{-1})v\xi}{qv\xi-u}-\frac{(q-q^{-1})v\xi}{qv\xi-(-u)}=
q^{-1}\big(\gle(qv\xi,u)-\gle(qv\xi,-u)\big)\quad\mbox{or}\\
&\quad=\frac{(q-q^{-1})q^{-1}u}{qv\xi-u}-\frac{(q-q^{-1})q^{-1}(-u)}{qv\xi-(-u)}=
q^{-1}\big(\gri(qv\xi,u)-\gri(qv\xi,-u)\big)\,.
\end{split}
\end{equation*}
Analogously, using \r{g-g} for both cases $\ell<j$ or $\ell>j$
the right hand side of \r{em18} is equal to
\begin{equation*}
\begin{split}
&\frac{2(q-q^{-1})v\xi u}{v^2\xi^2-u^2}\sk{1-
\frac{(q^2-1)u^2}{q^2v^2\xi^2-u^2}}=
\frac{2(q-q^{-1})q^2v\xi u}{q^2v^2\xi^2-u^2}=\\
&\quad=q\sk{\frac{(q-q^{-1})qv\xi}{qv\xi-u}-\frac{(q-q^{-1})qv\xi}{qv\xi-(-u)}}=
q\big(\gle(qv\xi,u)-\gle(qv\xi,-u)\big)\quad\mbox{or}\\
&\quad=q\sk{\frac{(q-q^{-1})u}{qv\xi-u}-\frac{(q-q^{-1})(-u)}{qv\xi-(-u)}}=
q\big(\gri(qv\xi,u)-\gri(qv\xi,-u)\big)\,.
\end{split}
\end{equation*}\qed

Let us define $\Rer^{(n)}(u,v)=\Id\ot\Id+\PPd^{(n)}(u,v)$, where operator $\PPd^{(n)}(u,v)$
is defined by \r{PP14}. Recall that $\RR$-matrix \r{tRalt} takes the form of sum of two
terms: $\RR^{(n)}(u,v)=\Rer^{(n)}(u,v)+\QQd^{(n)}(u,v)$,
where operator $\QQd^{(n)}(u,v)$ is defined by \r{QQ14}.

To prove equalities \r{e8}, we will consider only first equality \r{e8a} since
the proof of \r{e8b} is similar. To prove \r{e8a} we will prove separately
\begin{equation}\label{em11}
\begin{split}
&\RR^{(n)}_{12}(1,q)\RR^{(n)}_{34}(1,q)\RR^{(n)}_{14}(u,qv)\Rer^{(n)}_{13}(u,v)|i,1,j,1\>=\\
&\quad = \RR^{(n)}_{12}(1,q)\RR^{(n)}_{34}(1,q)\Rer^{(n-1)}_{13}(u,v)|i,1,j,1\>
\end{split}
\end{equation}
and
\begin{equation}\label{em12}
\begin{split}
&\RR^{(n)}_{12}(1,q)\RR^{(n)}_{34}(1,q)\RR^{(n)}_{14}(u,qv)\QQd^{(n)}_{13}(u,v)|i,1,j,1\>=\\
&\quad = \RR^{(n)}_{12}(1,q)\RR^{(n)}_{34}(1,q)\QQd^{(n-1)}_{13}(u,v)|i,1,j,1\>\,.
\end{split}
\end{equation}

Let $i^*=i$ for $i\not=n,n+1$ and $i^*=i'$ for $i=n,n+1$.
To prove \r{em11} we write the vector $\Rer^{(n)}_{13}(u,v)|i,1,j,1\>$
for $1<i,j<2n$ in the form
\begin{equation}\label{em13}
\Rer^{(n)}_{13}(u,v) |i,1,j,1\>=|i,1,j,1\>+\frac12\pf_{ij}(u,v)|j,1,i,1\>+\frac12\pf_{i^*j}(-u,v)|j^*,1,i^*,1\>\,.
\end{equation}

To calculate now $\RR^{(n)}_{14}(u,qv)\Rer^{(n)}_{13}(u,v) |i,1,j,1\>$
one can observe that only
$\Rer^{(n)}_{14}(u,qv)$ from $\RR^{(n)}_{14}(u,qv)$
acts to the vectors on the right hand side of
\r{em13} nontrivially. This is because the action of $\QQd^{(n)}_{14}(u,qv)$
is not vanishing only to the vectors of the form $|1',*,*,1\>$ which are absent
in \r{em13} since $1<j<2n$.
So the action of the operator $\RR^{(n)}_{14}(u,qv)$ to the all vectors
in the right hand sides of equalities \r{em13} add to these vectors the terms
which are linear combinations of the vectors $|1,1,k,l\>$ with $1 < k,l < 2n$.
But all the vectors of this form  are annihilated by the operator
$\RR^{(n)}_{12}(1,q)$ in the left hand side of \r{em11} due to the equality \r{em1}.
As result the action of the operator $\RR^{(n)}_{14}(u,qv)$ in the left hand
side of \r{em11} trivialises and we conclude the proof of the equality \r{em11} since
for $1<i,j<2n$ the vector $\Rer^{(n)}_{13}(u,v) |i,1,j,1\>$ coincides with the vector
$\Rer^{(n-1)}_{13}(u,v) |i,1,j,1\>$.

In order to prove equality \r{em12}, one has to consider the action of the matrix
$\QQd^{(n)}_{13}(u,v)$ to the vector $|i,1,j,1\>$ for $1<i,j<2n$. Using \r{QQ14}
one can check that $\QQd^{(n)}_{13}(u,v)|i,1,j,1\>=0$ for $i=n,n+1$, $j\not=n,n+1$
and for $i\not=n,n+1$, $j=n,n+1$. This means that \r{em12} is trivially satisfied
for these cases, and
to prove it, one has to consider the rest two cases
for $i,j\not=n,n+1$ and $i,j=n,n+1$.

Let us write explicitly the action of $\QQd^{(n)}_{13}(u,v)$ to the
vector $|i,1,j,1\>$ in these cases
\begin{itemize}
\begin{subequations}\label{em14}
\item $i\not=n,n+1\quad\mbox{and}\quad j\not=n,n+1$
\begin{equation}\label{em14a}
\begin{split}
&\QQd^{(n)}_{13}(u,v)|i,1,j,1\>=\frac{\delta_{ij'}}2\sum_{\ell=1\atop\ell\not=n,n+1}^{2n}
\Big(\qf_{\ell j}(u,v|\xi)+\qf_{\ell j}(-u,v|\xi)\Big)|\ell',1,\ell,1\>\,+\\
&\qquad+\frac{\delta_{ij'}}2\sum_{\ell=n}^{n+1}
\Big(\varphi^{(n)}_{\ell j}(u,v|\xi)|\ell',1,\ell,1\>+\varphi^{(n)}_{\ell j}(-u,v|\xi)|\ell',1,\ell',1\>\Big)\,,
\end{split}
\end{equation}
\item $i=n,n+1\quad\mbox{and}\quad j=n,n+1$
\begin{equation}\label{em14b}
\begin{split}
\QQd^{(n)}_{13}(u,v)|i,1,j,1\>&=\frac{1}2\sum_{\ell=1\atop\ell\not=n,n+1}^{2n}
\Big(\delta_{ij'}\varphi^{\prime(n)}_{\ell j}(u,v|\xi)+
\delta_{ij}\varphi^{\prime(n)}_{\ell j}(-u,v|\xi)\Big)|\ell',1,\ell,1\>\,+\\
&\ +\frac{\delta_{i,j'}}{4}\Big(\QQff(u,v|\xi)\Omega_{13}+
\QQdd(u,v|\xi)\Omega'_{13} - \omega_{13}^{(j)}\Big)+\\
&\ +\frac{\delta_{i,j}}{4}U_1 \Big(\QQff(-u,v|\xi)\Omega_{13}+
\QQdd(-u,v|\xi)\Omega_{13}' - \,\omega_{13}^{(j)}\Big)\,,
\end{split}
\end{equation}
\end{subequations}
\end{itemize}
where
\begin{equation}\label{TT15a}
\begin{split}
\varphi^{(n)}_{ij}(u,v|\xi)&=
\alpha_n^+(q)\ \qf_{ij}(u,v|\xi)+\alpha_n^-(q)\ \qf_{ij}(-u,v|\xi)\,,\\
\varphi^{\prime(n)}_{ij}(u,v|\xi)&=\alpha_n^+(q^{-1})\ \qf_{ij}(u,v|\xi)+
\alpha_n^-(q^{-1})\ \qf_{ij}(-u,v|\xi)\,,
\end{split}
\end{equation}
\begin{equation}
\alpha_n^\pm(q)=(1\pm q^{n-1})/2q^{1/2}, \end{equation}
and
\begin{equation*}
\Omega_{13}=\sum_{\ell=n}^{n+1}|\ell',1,\ell,1\>,\qquad
\Omega_{13}'=\sum_{\ell=n}^{n+1}|\ell',1,\ell',1\>, \qquad 
\omega_{13}^{(j)} = \alpha_q |j',1,j,1\> + \gamma_q \epsilon_j |j,1,j',1\>
\end{equation*}

In \r{em14a} and \r{em14b} we write explicitly dependence of the all the functions
 on $\xi=q^{-n+1}$. Our goal now is to demonstrate that the action of the
operator $\RR^{(n)}_{14}(u,qv)$ to the right hand sides of equalities
\r{em14a} and \r{em14b} and multiplication  from the left by the product of the operators
$\RR^{(n)}_{12}(1,q)\RR^{(n)}_{34}(1,q)$ results in the same formulas
with $\xi\to q\xi$ which is equivalent to the shift of $n\to n-1$ or
to the action of $\QQd^{(n-1)}_{13}(u,v)$ to the vectors $|i,1,j,1\>$ for all $1<i,j<2n$.

\begin{remark}\label{Remphi}
During the proof of equality \r{em12} there will appear functions \r{TT15a}
with shifted value $n$ and the same parameter $\xi$. This is a reason of notations
introduced by \r{TT15a} where we write separately dependence on $n$ through
coefficients $\alpha_n^\pm(q)$, $\alpha_n^\pm(q^{-1})$ and parameter
$\xi$ through the functions $\qf_{ij}(u,v|\xi)$.
\end{remark}

According to \r{D2Ralt}
the action of the operator $\RR^{(n)}_{14}(u,qv)$ to the vector $|\ell',1,\ell,1\>$ or
$|\ell',1,\ell',1\>$ consists of three parts
\begin{equation}\label{acqn}
\RR^{(n)}_{14}(u,qv)|\ell',1,\ell,1\>=|\ell',1,\ell,1\>+\PPd^{(n)}_{14}(u,qv)|\ell',1,\ell,1\>+
\QQd^{(n)}_{14}(u,qv)|\ell',1,\ell,1\>\,.
\end{equation}
Due to the structure of the operator $\PPd(u,v)$ \r{PP14} the  action
of $\PPd^{(n)}_{14}(u,qv)$ to the
vectors $|\ell',1,\ell,1\>$ and $|\ell',1,\ell',1\>$
will be always proportional to the vectors $|1,1,\ell,\ell'\>$ or $|1,1,\ell,\ell\>$.
But all such vectors are annihilated by the operator $\RR^{(n)}_{12}(1,q)$
in the left hand side of \r{em12}. It means that only the action of
$\QQd^{(n)}_{14}(u,qv)$ to the vectors $|\ell',1,\ell,1\>$ in \r{acqn} should be taken into account.

Proving \r{em12} we first consider the case $i,j\not=n,n+1$.
 The action of
the operator ${\QQd}^{(n)}_{13}(u,v)$ to the vectors $|i,1,j,1\>$
is given by \r{em14a}  and is non-vanishing
only if $i=j'$.
Next step is to act by ${\RR}^{(n)}_{14}(u,qv)$ onto right hand side of \r{em14a}.
So the vector
\begin{equation}\label{TT20}
{\RR}^{(n)}_{12}(1,q){\RR}^{(n)}_{34}(1,q){\RR}^{(n)}_{14}(u,qv)
{\QQd}^{(n)}_{13}(u,v)|i,1,j,1\>
\end{equation}
is equal to the sum
\begin{equation}\label{TT19}
\begin{split}
& \frac{\delta_{ij'}}2 {\RR}^{(n)}_{12}(1,q){\RR}^{(n)}_{34}(1,q)
\sum_{\ell=2\atop \ell\not=n,n+1}^{2n-1}
\Big(\qf_{\ell j}(u,v|\xi)+\qf_{\ell j}(-u,v|\xi)\Big)|\ell',1,\ell,1\>\ +\\
&+ \frac{\delta_{ij'}}2
{\RR}^{(n)}_{12}(1,q){\RR}^{(n)}_{34}(1,q)
\sum_{\ell=n}^{n+1}\Big(\varphi^{(n)}_{\ell j}(u,v|\xi)|\ell',1,\ell,1\>
+\varphi^{(n)}_{\ell j}(-u,v|\xi)|\ell',1,\ell',1\>\Big)
\end{split}
\end{equation}
and
\begin{equation}\label{TT18}
\begin{split}
&{\RR}^{(n)}_{12}(1,q){\RR}^{(n)}_{34}(1,q)
{\QQd}^{(n)}_{14}(u,qv){\QQd}^{(n)}_{13}(u,v)|i,1,j,1\>=\\
&\quad =
\frac{\delta_{ij'}}2\Big(\qf_{1j}(u,v|\xi)+\qf_{1j}(-u,v|\xi)\Big)
{\RR}^{(n)}_{12}(1,q){\RR}^{(n)}_{34}(1,q){\QQd}^{(n)}_{14}(u,qv)
|1',1,1,1\>\,.
\end{split}
\end{equation}
Note that in \r{TT19} the limits of summation  over $\ell$ are
changed. This is  because for $\ell=2n$ the vector $|1,1,2n,1\>$ is annihilated by the matrix
 ${\RR}^{(n)}_{12}(1,q)$ and for $\ell=1$ the vector $|2n,1,1,1\>$
 is annihilated by the matrix ${\RR}^{(n)}_{34}(1,q)$.
Also, the term with action $\PPd^{(n)}_{14}(u,qv)$ vanishes under the action of  ${\RR}^{(n)}_{12}(1,q)$.
 Using \r{em14a} again to calculate
 \begin{equation}\label{TT21}
\begin{split}
&{\QQd}^{(n)}_{14}(u,qv)|1',1,1,1\>=
 \frac{1}2 \sum_{\ell=1\atop \ell\not=n,n+1}^{2n}
\Big(\qf_{\ell 1}(u,qv|\xi)+\qf_{\ell 1}(-u,qv|\xi)\Big)|\ell',1,1,\ell\>\,+\\
&\qquad+ \frac{1}2 \sum_{\ell=n}^{n+1}\Big(\varphi^{(n)}_{\ell 1}(u,qv|\xi)|\ell',1,1,\ell\>
+\varphi^{(n)}_{\ell 1}(-u,qv|\xi)|\ell',1,1,\ell'\>\Big)
\end{split}
\end{equation}
and collecting the terms from \r{TT19} and \r{TT18} one concludes that the vector
\r{TT20} is equal to
\begin{equation}\label{TT22}
\frac{\delta_{ij'}}2 \RR^{(n)}_{12}(1,q)\RR^{(n)}_{34}(1,q)
\Big(A_j(u,v)+B_j(u,v)\Big)\,,
\end{equation}
where
\begin{equation}\label{TT23}
\begin{split}
A_j(u,v)&=\sum_{\ell=2\atop \ell\not=n,n+1}^{2n-1}
\Big((\qf_{\ell j}(u,v|\xi)+\qf_{\ell j}(-u,v|\xi))|\ell',1,\ell,1\>\,+\\
&\quad+ \frac12 (\qf_{1j}(u,v|\xi)+\qf_{1j}(-u,v|\xi))
(\qf_{\ell 1}(u,qv|\xi)+\qf_{\ell 1}(-u,qv|\xi))|\ell',1,1,\ell\>
\Big)
\end{split}
\end{equation}
and
\begin{equation}\label{TT24}
\begin{split}
&B_j(u,v)=\sum_{\ell=n}^{n+1}
\Big(\varphi^{(n)}_{\ell j}(u,v|\xi)|\ell',1,\ell,1\>+\varphi^{(n)}_{\ell j}(-u,v|\xi))|\ell',1,\ell',1\>\,+\\
&\quad+ \frac12 (\qf_{1j}(u,v|\xi)+\qf_{1j}(-u,v|\xi))
(\varphi^{(n)}_{\ell 1}(u,qv|\xi)|\ell',1,1,\ell\> +\varphi^{(n)}_{\ell 1}(-u,qv|\xi)|\ell',1,1,\ell'\>)
\Big)\,.
\end{split}
\end{equation}
Applying equality (cf. \r{em6})
\begin{equation*}
{\RR}^{(n)}_{34}(1,q)|\ell',1,1,\ell\>=
- q {\RR}^{(n)}_{34}(1,q)|\ell',1,\ell,1\>
\end{equation*}
for $1<\ell<n$ and $n+1<\ell<2n$ to the second line in \r{TT23} one concludes
using identity \r{em19} that
\begin{equation}\label{TT25}
A_j(u,v)=\sum_{\ell=2\atop \ell\not=n,n+1}^{2n-1}
\Big((\qf_{\ell j}(u,v|q\xi)+\qf_{\ell j}(-u,v|q\xi))|\ell',1,\ell,1\>\,.
\end{equation}

In order to transform vector $B_j(u,v)$ \r{TT24} we consider linear combination of the
vectors
\begin{equation}\label{TT26}
{\RR}^{(n)}_{34}(1,q)(\varphi^{(n)}_{\ell 1}(u,qv|\xi)|\ell',1,1,\ell\>
+\varphi^{(n)}_{\ell 1}(-u,qv|\xi)|\ell',1,1,\ell'\>)
\end{equation}
which can be rewritten in the form
\begin{equation*}
\begin{split}
&\frac12{\RR}^{(n)}_{34}(1,q)(\varphi^{(n)}_{\ell 1}(u,qv|\xi)+
\varphi^{(n)}_{\ell 1}(-u,qv|\xi))\Big(|\ell',1,1,\ell\> +|\ell',1,1,\ell'\>\Big)\,+\\
&\quad +\frac12{\RR}^{(n)}_{34}(1,q)(\varphi^{(n)}_{\ell 1}(u,qv|\xi)-
\varphi^{(n)}_{\ell 1}(-u,qv|\xi))\Big(|\ell',1,1,\ell\> -|\ell',1,1,\ell'\>\Big)
\end{split}
\end{equation*}
or (cf. \r{em99} and \r{em98})
\begin{equation*}
\begin{split}
&-\frac{q}2{\RR}^{(n)}_{34}(1,q)(\varphi^{(n)}_{\ell 1}(u,qv|\xi)+
\varphi^{(n)}_{\ell 1}(-u,qv|\xi))\Big(|\ell',1,\ell,1\> +|\ell',1,\ell',1\>\Big)\,+\\
&\quad -\frac12{\RR}^{(n)}_{34}(1,q)(\varphi^{(n)}_{\ell 1}(u,qv|\xi)-
\varphi^{(n)}_{\ell 1}(-u,qv|\xi))\Big(|\ell',1,\ell,1\> -|\ell',1,\ell',1\>\Big)\,.
\end{split}
\end{equation*}

Using explicit formulas of the functions $\varphi^{(n)}_{\ell m}(u,v|\xi)$
given by \r{TT15a} the vector \r{TT26} is equal to (see Remark~\ref{Remphi})
\begin{equation}\label{tt27}
\begin{split}
&{\RR}^{(n)}_{34}(1,q)(\varphi^{(n)}_{\ell 1}(u,qv|\xi)|\ell',1,1,\ell\>
+\varphi^{(n)}_{\ell 1}(-u,qv|\xi)|\ell',1,1,\ell'\>)=\\
&\quad=-q{\RR}^{(n)}_{34}(1,q)(\varphi^{(n-1)}_{\ell 1}(u,qv|\xi)|\ell',1,\ell,1\>
+\varphi^{(n-1)}_{\ell 1}(-u,qv|\xi)|\ell',1,\ell',1\>)
\end{split}
\end{equation}
and
\begin{equation}\label{TT27}
\begin{split}
&{\RR}^{(n)}_{34}(1,q) B_j(u,v)= {\RR}^{(n)}_{34}(1,q) \sum_{\ell=n}^{n+1}
\Big(\varphi^{(n)}_{\ell j}(u,v|\xi)|\ell',1,\ell,1\>+\varphi^{(n)}_{\ell j}(-u,v|\xi))|\ell',1,\ell',1\>\,-\\
&- \frac{q}2 \big(\qf_{1j}(u,v|\xi)+
\qf_{1j}(-u,v|\xi)\big)\big(\varphi^{(n-1)}_{\ell 1}(u,qv|\xi)|\ell',1,\ell,1\>+\varphi^{(n-1)}_{\ell 1}(-u,qv|\xi)|\ell',1,\ell',1\>\big)
\Big)\,.
\end{split}
\end{equation}
One can observe that coefficient functions in front of the vectors
$|\ell',1,\ell,1\>$ and $|\ell',1,\ell',1\>$ are related by the transformation
$u\to -u$. It is sufficient to consider the coefficient at one of this vector,
say $|\ell',1,\ell,1\>$.
This coefficient is equal to
\begin{equation}\label{TT30}
\varphi^{(n)}_{\ell j}(u,v|\xi)-\frac{q}2 (\qf_{1j}(u,v|\xi)+\qf_{1j}(-u,v|\xi))
\varphi^{(n-1)}_{\ell 1}(u,qv|\xi)\,.
\end{equation}
Using identities \r{em19} and \r{em20}
proved in Lemma~\ref{iden} one can check that combination
of the functions \r{TT30} is equal to
\begin{equation}\label{TT31}
\varphi^{(n-1)}_{\ell j}(u,v|q\xi)=\alpha^+_{n-1}(q)\qf_{\ell j}(u,v|q\xi)+
\alpha^-_{n-1}(q)\qf_{\ell j}(-u,v|q\xi)
\end{equation}
and the vector $B_j(u,v)$ is
\begin{equation}\label{TT28}
{\RR}^{(n)}_{34}(1,q)  B_j(u,v)={\RR}^{(n)}_{34}(1,q) \sum_{\ell=n}^{n+1}
\Big(\varphi^{(n-1)}_{\ell j}(u,v|q\xi)|\ell',1,\ell,1\>+
\varphi^{(n-1)}_{\ell j}(-u,v|q\xi))|\ell',1,\ell',1\>\Big)\,.
\end{equation}
Taking into account that product of $\xi$ and $q$ means effectively the change $n\to n-1$, summing up \r{TT25} and \r{TT28} and comparing with \r{em14a} one can see that the vector \r{TT20} is equal to
\begin{equation*}
{\RR}^{(n)}_{12}(1,q){\RR}^{(n)}_{34}(1,q){\QQd}^{(n-1)}_{13}(u,v)|i,1,j,1\>
\end{equation*}
and so equality \r{em12} is proved for $i,j\not=n,n+1$.

For $i,j=n,n+1$ the vector $\QQd^{(n)}_{13}(u,v)|i,1,j,1\>$
is given by  \r{em14b}.
Let us calculate first the coefficients at the vectors
 $|\ell',1,\ell,1\>$  for $1\leq \ell<n$ and $n+1<\ell\leq 2n$
 in
\begin{equation}\label{lhs522}
{\RR}^{(n)}_{12}(1,q){\RR}^{(n)}_{34}(1,q){\RR}^{(n)}_{14}(u,qv)
\QQd^{(n)}_{13}(u,v)|i,1,j,1\>\,.
\end{equation}
These vectors proportional to $\delta_{ij'}$ are given by the equality
\begin{equation}\label{TT44}
\begin{split}
&\frac12{\RR}^{(n)}_{12}(1,q){\RR}^{(n)}_{34}(1,q)
\sum_{\ell=2\atop \ell\not=n,n+1}^{2n-1}\varphi^{\prime(n)}_{\ell j}(u,v|\xi)|\ell',1,\ell,1\>\,+\\
&-\frac{q}4{\RR}^{(n)}_{12}(1,q){\RR}^{(n)}_{34}(1,q)
\varphi^{\prime(n)}_{1 j}(u,v|\xi)\sum_{\ell=2\atop \ell\not=n,n+1}^{2n-1}
\big(\qf_{\ell 1}(u,qv|\xi)+\qf_{\ell 1}(-u,qv|\xi)\big)
|\ell',1,\ell,1\>,
\end{split}
\end{equation}
where we have used the action \r{TT21} and
${\RR}^{(n)}_{34}(1,q)|\ell',1,1,\ell\>=-q{\RR}^{(n)}_{34}(1,q)|\ell',1,\ell,1\>$.
Gathering the coefficients at the same vectors in \r{TT44} and
using the explicit form of the functions $\varphi^{\prime(n)}_{\ell j}(u,v)$ one obtains
from the rational identities \r{em19} and \r{em18} that
\begin{equation}\label{TT45}
\varphi^{\prime(n)}_{\ell j}(u,v)-\frac{q}{2}\varphi^{\prime(n)}_{1 j}(u,v)
\big(\qf_{\ell 1}(u,qv)+\qf_{\ell 1}(-u,v)\big)=\varphi^{\prime(n-1)}_{\ell j}(u,v)\,.
\end{equation}
Analogous terms coming from \r{lhs522} which are proportional to $\delta_{ij}$
have coefficients
$\varphi^{\prime(n-1)}_{\ell j}(-u,v)$ and this proves that the coefficients at
the vectors $|\ell',1,\ell,1\>$  for $\ell\not=n,n+1$
in the  right and left
hand sides in the equality \r{em12} coincide.

The vectors which are proportional to $\delta_{ij'}/4$ in the vector
\r{lhs522} are
\begin{equation}\label{TT50}
\begin{split}
&\Big(\QQff(u,v|\xi)-q\varphi^{\prime(n)}_{1j}(u,v|\xi)\varphi^{(n-1)}_{n1}(u,qv|\xi)\Big)
\Omega_{13}\ +\\
&\qquad+\Big(\QQdd(u,v|\xi)-q\varphi^{\prime(n)}_{1j}(u,v|\xi)
\varphi^{(n-1)}_{n1}(-u,qv|\xi)\Big)\Omega'_{13} - \omega_{13}^{(j)}=\\
&\qquad= \QQff(u,v|q\xi)\Omega_{13}+
\QQdd(u,v|q\xi)\Omega'_{13} - \omega_{13}^{(j)}
\end{split}
\end{equation}
The matrix $\QQd^{(n)}_{14}(u,qv)$ acts nontrivially only to the single
vector $|1',1,1,1\>$ from the first line of \r{em14b}. Then to obtain left hand side of
\r{TT50}  one has to use equalities \r{TT21} and \r{tt27}. Equality in
\r{TT50} follows from the explicit form of the functions  \r{TT15a} and
\r{QQff}, which result in the following relations for $j = n, n+1$
\begin{equation}\label{TT52}
\QQff(u,v|\xi)-q\varphi^{\prime(n)}_{1j}(u,v|\xi)\varphi^{(n-1)}_{n1}(u,qv|\xi)=
\QQff(u,v|q\xi)
\end{equation}
and
\begin{equation}\label{TT53}
\QQdd(u,v|\xi)-q\varphi^{\prime(n)}_{1j}(u,v|\xi)
\varphi^{(n-1)}_{n1}(-u,qv|\xi)=\QQdd(u,v|q\xi)\,.
\end{equation}

Analogously, the vectors which are proportional to $\delta_{ij}/4$ in the vector
\r{lhs522} are
\begin{equation}\label{TT51}
\begin{split}
&\Big(\QQff(-u,v|\xi)-q\varphi^{\prime(n)}_{1j}(-u,v|\xi)\varphi^{(n-1)}_{n1}(-u,qv|\xi)\Big)
\Omega'_{13}\ +\\
&\qquad+\Big(\QQdd(-u,v|\xi)-q\varphi^{\prime(n)}_{1j}(-u,v|\xi)
\varphi^{(n-1)}_{n1}(u,qv|\xi)\Big)\Omega_{13} - U_1\,\omega_{13}^{(j)}=\\
&\qquad=\QQff(-u,v|q\xi)\Omega_{13}'+
\QQdd(-u,v|q\xi)\Omega_{13} - U_1\,\omega_{13}^{(j)}\,.
\end{split}
\end{equation}

This proves that equality \r{e8a} is valid for all values of $1<i,j<2n$. \qed

\section{The composed currents and the Serre relations}\label{ApE}

In this appendix, we prove the proposition.
\begin{prop}
The Serre relations \r{SerreDn} are equivalent to the commutation
relations between simple and composed currents ($a=n,n+1$)
\begin{equation}\label{SRD}
\begin{split}
F_{n-2}(v)\ \Fn_{a,n-2}(u)&=f(v^2,u^2)\ \Fn_{a,n-1}(u)\ F_{n-2}(v)\,,\\
f(v^2,u^2)\ E_{n-2}(v)\ \En_{n-1,a}(u)&=\En_{n-1,a}(u)\ E_{n-2}(v)\,,\\
\Fn_{n+2,n-2}(u)\ F_{n-1}(v)&=f(u,-v)\ F_{n-1}(v)\ \Fn_{n+2,n-2}(u)\,,\\
f(u,-v)\ \En_{n-2,n+2}(u)\ E_{n-1}(v)&=E_{n-1}(v)\ \En_{n-2,n+2}(u)\,.
\end{split}
\end{equation}
\end{prop}
This proposition may be used to calculate the full algebra of the composed currents
\r{comp-cuF} and \r{comp-cuE}, which is necessary in the application to the
quantum integrable models associated with the algebra $U_q(D^{(2)}_{n})$.

{\it Proof.}  We will prove only the equivalence of the first and third Serre relations in \r{SerreDn}
to the first and third commutation relations in \r{SRD}  since other
equivalences can be proved analogously.
To simplify the consideration, we will consider them in the most
simple non-trivial case of the algebra $U_q(D^{(2)}_3)$ or $n=3$. We will use an approach
introduced in the paper \cite{DingKhoroshkin2000Weylgroup} and start from the
first Serre relation in \r{SerreDn} rewriting it in the form
\begin{equation}\label{AE1}
\mathop{\rm Sym}_{v_1,v_2}\Big(
F_1(v_1)F_1(v_2)F_2(u)-(q+q^{-1})F_1(v_1)F_2(u)F_1(v_2)+F_2(u)F_1(v_1)F_1(v_2)\Big)=0\,.
\end{equation}

To transform this relation, we will need the commutation relations between the currents
$F_1(v)$ and $F_2(u)$ (see \r{kEFD}) in the form
\begin{equation}\label{AE2}
F_1(v)F_2(u)=f(u^2,v^2)F_2(u)F_1(v)+c\Big(\delta(v,u)\Fn_{3,1}(u)+\delta(v,-u)\Fn_{4,1}(-u)\Big)\,,
\end{equation}
where $c=(q^{-1}-q)/2$ and composed currents $\Fn_{3,1}(u)=F_2(u)F_1(u)$
and $\Fn_{4,1}(u)=F_2(-u)F_1(u)$ are given by \r{comp-cuF} for the case $n=3$.

In the category of the highest weight representations, the matrix elements of the
 product of the currents
$\<w^*,F_{1}(v)F_{2}(u)w\>$ is a series which converges to some rational function
in the domain $|v|\gg|u|$. Here, $w$ is a vector from the highest weight representation
$W$ of $U_q(\Dn{3})$ and $w^*$ is a  vector from the dual representation $W^*$. Then
equality \r{AE2} between formal series can be regarded as an analytical
continuation from the domain $|v|\gg|u|$ to the domain $|v|\ll|u|$ where
the first term in the right hand side of \r{AE2} is well defined.
In this category of the representations of the quantum loop algebra $U_q(\Dn{3})$
the product of the currents $F_1(v)F_2(u)$ will have the simple poles at the
points $v=\pm u$ and the coefficients in front of the $\delta$-functions in \r{AE2}
are the residues of this product in the corresponding point. One can calculate these
residues using the commutation relations between currents
\begin{equation}\label{AE3}
\begin{split}
\mathop{\rm res}_{v=u}\ F_1(v)F_2(u)\frac{dv}{v}&=
\left.\frac{v-u}{v}\ F_1(v)F_2(u)\right|_{v=u}=
\left.\frac{1}{v}\frac{q^{-1}v^2-qu^2}{v+u}\ F_2(u)F_1(v)\right|_{v=u}=\\
&=c\ F_2(u)F_1(u)=c\ \Fn_{3,1}(u)\,.
\end{split}
\end{equation}
Analogously one can calculate
\begin{equation}\label{AE4}
\mathop{\rm res}_{v=-u}\ F_1(v)F_2(u)\frac{dv}{v}=
\left.\frac{v+u}{v}\ F_1(v)F_2(u)\right|_{v=-u}=c\ F_2(u)F_1(-u)=c\ \Fn_{4,1}(-u)\,.
\end{equation}

Substituting \r{AE2} to \r{AE1} one gets two types of terms. There are terms without
$\delta$-functions and the terms with them. The term without $\delta$-functions
will be  proportional to the product of the currents $F_2(u)F_1(v_1)F_1(v_2)$
with a coefficient of proportionality
\begin{equation*}
\Big(\Big(f(u^2,v_2^2)-(q+q^{-1})\Big)f(u^2,v^2_1)+1\Big)=(g(u^2,v^2_2)
-g(u^2,v^2_1))\ f(v_2^2,v_1^2)\,.
\end{equation*}
Due to the commutation relation between the currents $F_{n-2}(v_{1})$ and $F_{n-2}(v_{2})$ from \r{kEFD} this terms will disappear
after symmetrization over $v_1,v_2$. Equating the coefficients in front of the
$\delta$-functions on gets that the Serre relation \r{AE1} is equivalent
to the commutation relations
\begin{equation}\label{AE5}
F_1(v)\ \Fn_{a,1}(u)=f(v^2,u^2)\ \Fn_{a,1}(u)\ F_1(v),\quad a=3,4
\end{equation}
proving the first line in \r{SRD}.

The next step is to calculate the commutation relations of the composed current
$\Fn_{3,1}(v)$ and $\Fn_{4,1}(v)$ with the currents $F_2(u)$. Using \r{AE2}
one gets
\begin{equation}\label{AE6}
\begin{split}
\Fn_{3,1}(v)\ F_2(u)&=F_2(v)F_1(v)\ F_2(u)=\\
&=F_2(v)\Big(
f(u^2,v^2)F_2(u)F_1(v)+c\Big(\delta(v,u)\Fn_{3,1}(u)+\delta(v,-u)\Fn_{4,1}(-u)\Big)\Big)=\\
&=\fgo(u,v)\fgo(-u,v)\ F_2(v)F_2(u)F_1(v)+c\delta(v,-u)\ F_2(v) F_2(-v)F_1(v)=\\
&= \fgo(v,u)\fgo(-u,v)\ F_2(u)\Fn_{3,1}(v)-c\delta(v,-u)\ \Fn_{5,1}(v)\,,
\end{split}
\end{equation}
where we have used the fact that $F_2(v)^2=0$, commutativity of the currents $F_2(u)$
and $F_2(-u)$ and the commutation relation between the currents $F_2(v)$, $F_2(u)$.
Analogously, one can calculate that
\begin{equation}\label{AE7}
\Fn_{4,1}(v)\ F_2(u)= \fgo(u,v)\fgo(-v,u)\ F_2(u)\Fn_{4,1}(v)-c\delta(v,u)\ \Fn_{5,1}(v)\,.
\end{equation}
In \r{AE6} and \r{AE7} $\Fn_{5,1}(v)=-F_2(-v)F_2(v)F_1(v)$ is a composed current
given by the last line in \r{comp-cuF} for $n=3$.

We will use the commutation relations \r{AE2}, \r{AE6} and \r{AE7} to analyze
the third Serre relation in \r{SerreDn} which can be rewritten in the form
\begin{equation}\label{AE8}
\begin{split}
&\mathop{\rm Sym}_{v_1,v_2,v_3}\Big(
F_1(u)F_2(v_1)F_2(v_2)F_2(v_3)-[3]_q\ F_2(v_1)F_1(u)F_2(v_2)F_2(v_3)+\\
&\quad +[3]_q\ F_2(v_1)F_2(v_2)F_1(u)F_2(v_3)-F_2(v_1)F_2(v_2)F_2(v_3)F_1(u)\Big)=0\,,
\end{split}
\end{equation}
where
\begin{equation*}
[3]_q=q+1+q^{-1}\,.
\end{equation*}
Let us denote the combination of the composed currents and the $\delta$-functions
in \r{AE2} as
\begin{equation}\label{AE9}
\mathbb{F}(u,v)=c\Big(\delta(u,v)\ \Fn_{3,1}(u)+\delta(u,-v)\ \Fn_{4,1}(u)\Big)\,.
\end{equation}
Using the commutation relations \r{AE6} and \r{AE7} one can calculate
\begin{equation}\label{AE10}
\mathbb{F}(u,v_1)\ F_2(v_2)=\fgo(v_1,v_2)\fgo(-v_2,v_1)\ F_2(v_2)\ \mathbb{F}(u,v_1)+
\mathbb{\tilde{F}}(u,v_1,v_2)\,,
\end{equation}
where
\begin{equation}\label{AE11}
\mathbb{\tilde{F}}(u,v_1,v_2)=-2c^2\ \delta(v_1,-v_2)\ \delta(u^2,v_1^2)\ \Fn_{5,1}(u)\,.
\end{equation}
Here we have
introduced another combination of the product of two $\delta$-functions and the
composed current $\Fn_{5,1}(u)$. Using in \r{AE8} the commutation
relations \r{AE2} and \r{AE10} to move all the currents depending on
the spectral parameter $u$ to the right, we will obtain three types
of terms: the terms without the $\delta$-functions, the terms which
are proportional to the single  $\delta$-function and finally the terms
which are proportional to the product of two $\delta$-functions.
Due to the linear independence of the $\delta$-functions these terms
should vanish separately due to the Serre relation \r{AE8}.

We start considering the terms of the third type. Due to the identity
\begin{equation*}
[3]_q-\fgo(u,v)\fgo(-v,u)-f(v^2,u^2)=f(u,-v)
\end{equation*}
these terms are symmetrization over the set of the spectral parameters
$v_1,v_2,v_3$ of the following combination
\begin{equation*}
-2c^2\ \delta(v_1,-v_2)\ \delta(u^2,v_3^2)
\Big(\Fn_{5,1}(u)\ F_2(v_3)-f(u,-v_3)\ F_2(v_3)\ \Fn_{5,1}(u)\Big)\,.
\end{equation*}
Vanishing of these terms leads to the commutation relation
\begin{equation}\label{AE12}
\Fn_{5,1}(u)\ F_2(v)=f(u,-v)\ F_2(v)\ \Fn_{5,1}(u)\,.
\end{equation}

The terms of the second type are symmetrization of the following combination
\begin{equation*}
X(v_1,v_2,v_3)\ \fgo(v_2,v_1)\ F_2(v_1)F_2(v_2)\ \mathbb{F}(u,v_3)\,,
\end{equation*}
where the rational function
\begin{equation*}
X(v_1,v_2,v_3)=2q^{-1/2}c\, v_3\ (v_1-v_2)\frac{(1+q)v_3(v_1+v_2)-v_1v_2-qv_3^2}
{(v_1^2-v_3^2)(v_2^2-v^2_3)}
\end{equation*}
is anti-symmetric at $v_1\leftrightarrow v_2$. Because of this fact, the terms
of the second type vanish due to the commutation relations between
currents $F_2(v_1)$ and $F_2(v_2)$.

Finally, the first type term i symmetrization of the the combination
\begin{equation}\label{AE14}
\frac{2cu^2\ Y(u,v_1,v_2,v_3)}{q^2(u^2-v^2_1)(u^2-v^2_2)(u^2-v^2_3)}\
F_2(v_1)F_2(v_2)F_2(v_3)F_1(u)\,,
\end{equation}
where
\begin{equation*}
Y(u,v_1,v_2,v_3)=(q^2v^2_2-v^2_1)(q^2v_3-u^2)+q(q^2v^2_3-v^2_2)(u^2-v_1^2)\,.
\end{equation*}
One may verify that symmetrization of \r{AE14} over all spectral parameters
$v_1,v_2,v_3$ vanish due to the commutation relations between currents $F_2(v_i)$.
This proves that the Serre relation \r{AE8} is equivalent to the commutation relation
\r{AE12}.\qed

\bibliographystyle{JHEP}

\end{document}